\documentclass[12pt,a4paper]{amsart}
\usepackage[T1]{fontenc}
\usepackage[utf8]{inputenc}
\usepackage{amsmath}
\usepackage{amsfonts}
\usepackage{amssymb}
\usepackage{setspace}
\usepackage{amsthm}
\usepackage{booktabs}
\usepackage{rotating}
\usepackage[a4paper, left=2.5cm, right=2.5cm, top=2.8cm, bottom=2.8cm]{geometry}
\usepackage{float}
\usepackage{amsmath}
\usepackage{arydshln}
\usepackage{graphicx}
\usepackage[usenames, dvipsnames]{color}
 \usepackage{tabu}
\usepackage{subfigure,bm}
\usepackage{natbib}
\usepackage{todonotes}
\usepackage{commath}
\usepackage{amssymb}
\usepackage{ mathrsfs }
\usepackage{amsaddr} 
\usepackage{hyperref} 
\usepackage{url} 
\usepackage{comment} 

\newcommand\Ex{\text{\sf E}}
\newcommand\Prob{\text{\sf P}}

\newcommand\Cov{\text{\sf Cov}}

\newtheorem{theorem}{Theorem}
\newtheorem{proposition}{Proposition}
\newtheorem{lemma}{Lemma}
\newtheorem{corollary}{Corollary}

\theoremstyle{remark}
\newtheorem{example}{Example}
\newtheorem{remark}{Remark}

\begin{document}
\title[IIA through the Slepian model]{The Slepian model based independent interval approximation of persistency and zero-level excursion distributions}

\author[H.Bengtsson and K. Podg\'orski]{Henrik Bengtsson and Krzysztof Podg\'orski} \address{Department of Statistics,
Lund University, \\}
\email{Henrik.Bengtsson@stat.lu.se, \href{https://orcid.org/0000-0002-9280-4430}{ \includegraphics[height=2.2mm]{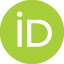} 0000-0002-9280-4430}}  
\email{Krzysztof.Podgorski@stat.lu.se, \href{https://orcid.org/0000-0003-0043-1532}{ \includegraphics[height=2.2mm]{ORCIDiD64.png} 0000-0003-0043-1532}}  
	\date{January 27, 2025}

\begin{abstract}
In physics and engineering literature, the distribution of the excursion time of a stationary Gaussian process has been approximated through a method based on a stationary switch process with independently distributed switching times. The approach matches the covariance of the clipped Gaussian process with that of the stationary switch process. By expressing the switching time distribution as a function of the covariance, the so-called independent interval approximation (IIA) is obtained for the excursions of Gaussian processes.  This approach has successfully approximated the persistency coefficient for many vital processes in physics but left an unanswered question about when such an approach leads to a mathematically meaningful and proper excursion distribution. Here, we propose an alternative approximation: the Slepian-based IIA. This approach matches the expected values of the clipped Slepian process and the corresponding switch process initiated at the origin. It is shown that these two approaches are equivalent, and thus, the original question of the mathematical validity of the IIA method can be rephrased using the Slepian model setup. We show that this approach leads to valid approximations of the excursion distribution for a large subclass of the Gaussian processes with monotonic covariance. Within this class, the approximated excursion time distribution has a stochastic representation that connects directly to the covariance of the underlying Gaussian process. This representation is then used to approximate the persistency coefficient for several important processes to illustrate the Slepian-based IIA approach. Lastly, we argue that the ordinary IIA is ill-suited in certain situations, such as for Gaussian processes with a non-monotonic covariance. 
\end{abstract}

\maketitle
\section{Introduction}
\noindent
The distribution of excursion times for a stationary Gaussian process constitutes a long-standing and challenging problem in physics and applied probability theory. In its most generic form, the problem can be formulated as finding statistical properties of the excursion set  $\mathcal{E}_u=\{t \in \mathbb R: X(t)>u\}$, where $X(t)$, $t\in \mathbb R$, is a stochastic process. This excursion is the exceedance above the level $u$, and thus, one often uses terms the {\it  $u$-excursion sets} and the {\it excursion distributions} when probabilistic properties of such sets are sought.

While the idea of an excursion distribution is intuitive, the solution to the problem of fully characterizing it has so far been elusive despite considerable effort. This effort has, however, resulted in several significant results, such as the famous Rice formula. Due to the importance of excursion distributions in many fields and the general complexity of the problem, several approximation methods have been developed. Two methods have been a frequent topic in physics papers. The first is based on Rice series expansions, which use the moments of the distribution of the number of crossings. The second method is the independent interval approximation (IIA).

The basis for the IIA method is the so-called clipped process. This process is obtained by computing the sign of a stochastic process, and \cite{McFadden1956} noted that the interval lengths of the clipped process correspond to the zero-level excursions of the original process. However, due to the underlying dependency of the original process, the intervals of the clipped process will have a non-trivial dependency structure. Therefore, in the IIA framework, the core assumption is that the dependence between successive intervals is sufficiently small so that they can be treated as independent. Subsequently, the covariance of the clipped process can be matched with that of a stationary binary renewal process, for which the relation between covariance and interval distributions is known. Using this relation, an approximation of the excursion distribution is obtained as a function of the covariance of the clipped process. While this approach leads to an approximation, the validity of such a method needs to be investigated.

The validity of the IIA approach was further discussed in \cite{McFadden1957AMS, McFadden1958}, which also considered extending the approach to allow for Markovian types of dependency between intervals. 
This extension led to a study by \cite{Rainal1962, Rainal1963}, which empirically examined both the IIA and the Markovian approximation. This study concluded that a Markovian type of dependency structure might not be suitable for approximating some Gaussian processes. Additionally, in \cite{LonguetHiggins1962}, it was shown that the IIA is never exactly valid for Gaussian processes. These investigations imply that the independence assumptions have to be used with caution. While these problems were highlighted, it did not mean that the approximation could be adequate for many  important models.

The popularity of the IIA has seen a resurgence during the last three decades. This is due to a greater interest in the tail of the excursion distribution in several areas such as optics, statistical physics, and more, \citep{Brainina2013}. To characterize the tail of this distribution, emphasis has been placed on the so-called {\sl persistency coefficient}. This coefficient characterizes the tail of the excursion distribution and thus describes the probability of large excursions. A comprehensive overview of persistency coefficients and how they can be approximated from the IIA framework can be found in \cite{BrayMS}.

Persistency coefficients are generally process-specific and often very difficult to derive analytically. Hence, they are only known explicitly for a handful of processes. The persistency coefficient for the Ornstein-Uhlenbeck process was shown to be $\theta=0.5$ by \cite{Bachelier} using the explicit distribution of the Brownian motion excursion above $u<0$, which is the Lévy distribution with the index of stability $0.5$. Another process for which the persistency coefficient is known is the random acceleration process. It was shown by \cite{Sinai:1992aa} that the persistency coefficient for this process is $\theta=0.25$. For fractional Brownian motions, which in its stationary version correspond to fractional Ornstein-Uhlenbeck processes, \cite{Molchan} showed that the persistency coefficient for this process is $\theta=1-H$, where $H$ is the Hurst parameter. Another important class of processes in statistical physics is diffusion processes. It was only in $2018$ that \cite{PoplavskyiSchehr2018} derived the persistency coefficient for the two-dimensional diffusion process $(\theta=0.1875)$. For this process, the IIA framework was used several times to approximate the persistency coefficient before it was derived. Estimates have been obtained by, among others \cite{Sire2007, Sire2008, PhysRevLett.99.060603, MajumdarSBC}, and they are relatively close to the persistency coefficient derived by \cite{PoplavskyiSchehr2018}. The importance of the diffusion process and the recency of the result by \cite{PoplavskyiSchehr2018} shows that there is still a need for adequate approximation methods. Even for the zero-level excursions.

For the cases with analytically derived persistency coefficients, the IIA generally seems to produce reasonable approximations. However, several important questions still remain that have not been addressed so far. Since the approximated excursion distribution is given as a function of the covariance of the clipped process in the Laplace domain, this function needs to be completely monotone for it to correspond to a valid probability distribution. It is not apparent that the class of covariance functions of the clipped Gaussian process will always correspond to a probability distribution on a positive real line. This also raises questions about the validity of using the IIA to approximate the persistency coefficient.

The main idea of this paper is to base the IIA on the Slepian process, which describes the behavior at the crossing instance of the process we are interested in. The reasons for using the clipped Slepian process are twofold; firstly, by matching the clipped Slepian process to a non-stationary binary process, the so-called inspection paradox can be avoided since all intervals of the non-stationary binary process will have the same interval distribution. Secondly, in a recent paper by \cite{Bengtsson2}, conditions were derived for when the interval distribution of a non-stationary binary can be uniquely obtained from the expected value function of this binary process. The latter result will be used to show that the approximated excursion distribution of the Slepian-based IIA is indeed a valid probability distribution for a large class of processes. This has direct implications for the ordinary IIA since it will later be shown that the ordinary IIA and the Slepian-based IIA are equivalent for zero-level excursions. Before introducing the paper's notation and outline of the paper, it should be noted that only processes with finite crossing intensities are treated unless otherwise stated. This is the same general assumption on crossing intensity that \cite{BrayMS} gives for the ordinary IIA.

The structure of the paper is as follows: we start with an overview of the ordinary IIA and introduce the switch process in Section \ref{sec:IIA}. Then, the Slepian-based IIA is the topic of Section \ref{sec:slepianIIA}. These sections focus on the approximate excursion distribution while Section \ref{sec:Montecarlo} focuses on the persistency coefficient. In Section \ref{sec:examples}, these methods are used to obtain approximated persistency coefficients for several stochastic processes. Due to the idea of matching and clipping different processes, the notation used in the paper is summarized in Table~\ref{tab:IIA-principles}.
\begin{table}[t!]
    \centering
    \begin{tabular}{lcc}\toprule
      & \sc Clipped process  & \sc Binary process \\ 
      \midrule[1pt]
    & \multicolumn{2}{c}{\bf Stationary}
     \\
    \midrule
    Process &  $D_{cl}(t)\stackrel{ \mbox{\tiny def}}{=}{\rm sgn}(X(t))$   &
     $D_s(t)$ \\
     \midrule
     Covariance  &  $R_{cl}(t)\stackrel{\mbox{\tiny def}}{=}\Cov\left(D_{cl}(t),D_{cl}(0)\right)$ &
     $R(t)\stackrel{\mbox{\tiny def}}{=}\Cov\left(D_s(t),D_s(0)\right)$\\
    \midrule[1pt]
    & \multicolumn{2}{c}{\bf Slepian}
     \\
    \midrule
    Process  & $D_0(t)\stackrel{\mbox{\tiny def}}{=}{\rm sgn}(X_0(t))$  & 
     $D(t)$ 
     \\
     \midrule
    Expectation &  $E_0(t)\stackrel{\mbox{\tiny def}}{=}\Ex\left(D_0(t)\right)$ &  $E(t)\stackrel{\mbox{\tiny def}}{=}\Ex\left( D(t)\right) $\\ 
     \bottomrule \vspace{1mm}
    \end{tabular}
    \caption{Notation for the characteristics that are used for the IIA matching clipped processes with switch processes in the stationary and attached at zero-crossing cases.}
    \label{tab:IIA-principles}
\end{table}

\section{The independent interval approximation} \label{sec:IIA}
\noindent
The instants of zero-level crossings of a smooth process $X(t)$ form an ordered sequence of points $S_i$, $i\in \mathbb Z$, leading to two sequences of interlaced intervals of the lengths: $T_i^+$, $i\in \mathbb Z$, for the excursions above the zero-level and $T_i^-$, $i\in \mathbb Z$, for the excursions below it. If a new process is constructed from this process such that it takes the value one when $X(t)\geqslant 0 $ and minus one when $X(t)<0$, we obtain the clipped version of $X(t)$. The clipped version can also be obtained directly from $X(t)$ by computing
\begin{align*}
   D_{cl}(t)={\rm sign}(X(t)). 
\end{align*} 
It is intuitive that the distributions of the intervals that $D_{cl}$ takes the value one or minus one encapsulate the excursion behavior of $X(t)$. This has been known for a long time. See, for example, \cite{McFadden1956}.

\begin{figure}[t]
\includegraphics[width=0.80\textwidth]{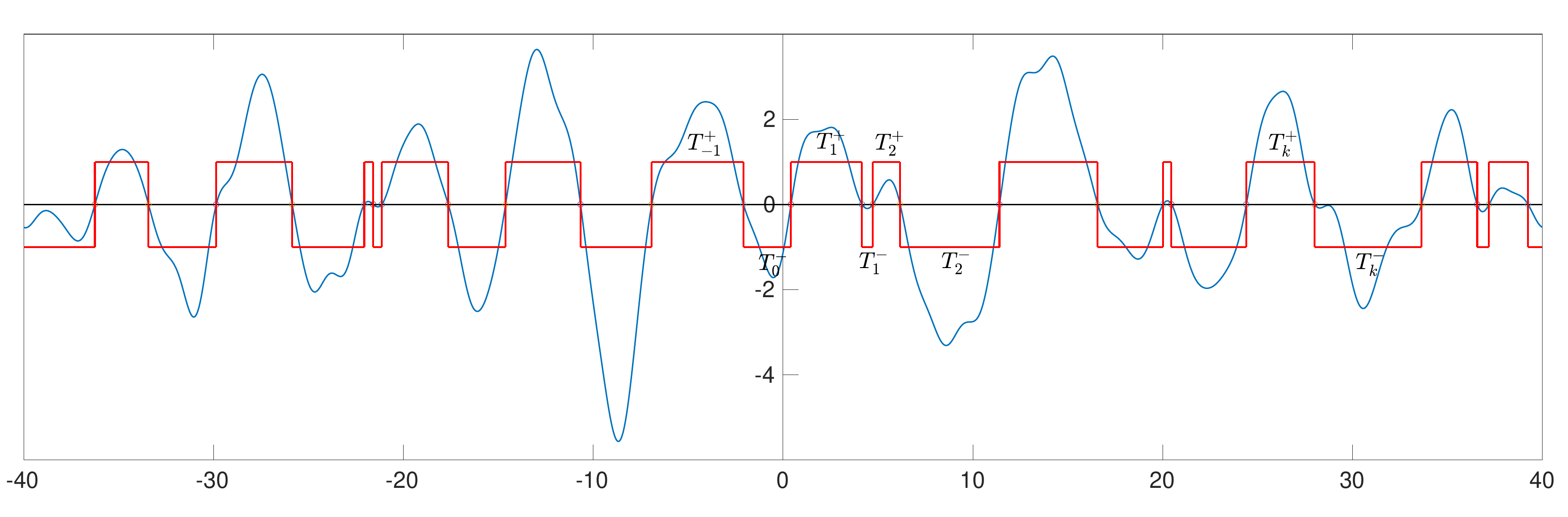}\\
\caption{\small The excursion intervals of a Gaussian process $X(t)$ together with the corresponding clipped process.}
\label{swcl}
\end{figure}

While the notion of clipping a process $X(t)$ might be straightforward, we illustrate it with an example, in which we use a double peak spectral model for sea state data that accounts for dependence between the two sea systems: swell and wind. 
The model was developed based on measured spectra for Norwegian waters (Haltenbanken and Statfjord); see \cite{Torsethaugen}. 
\begin{figure}[t]
\rule{-1.8cm}{0cm}
{\includegraphics[width=0.55\textwidth]{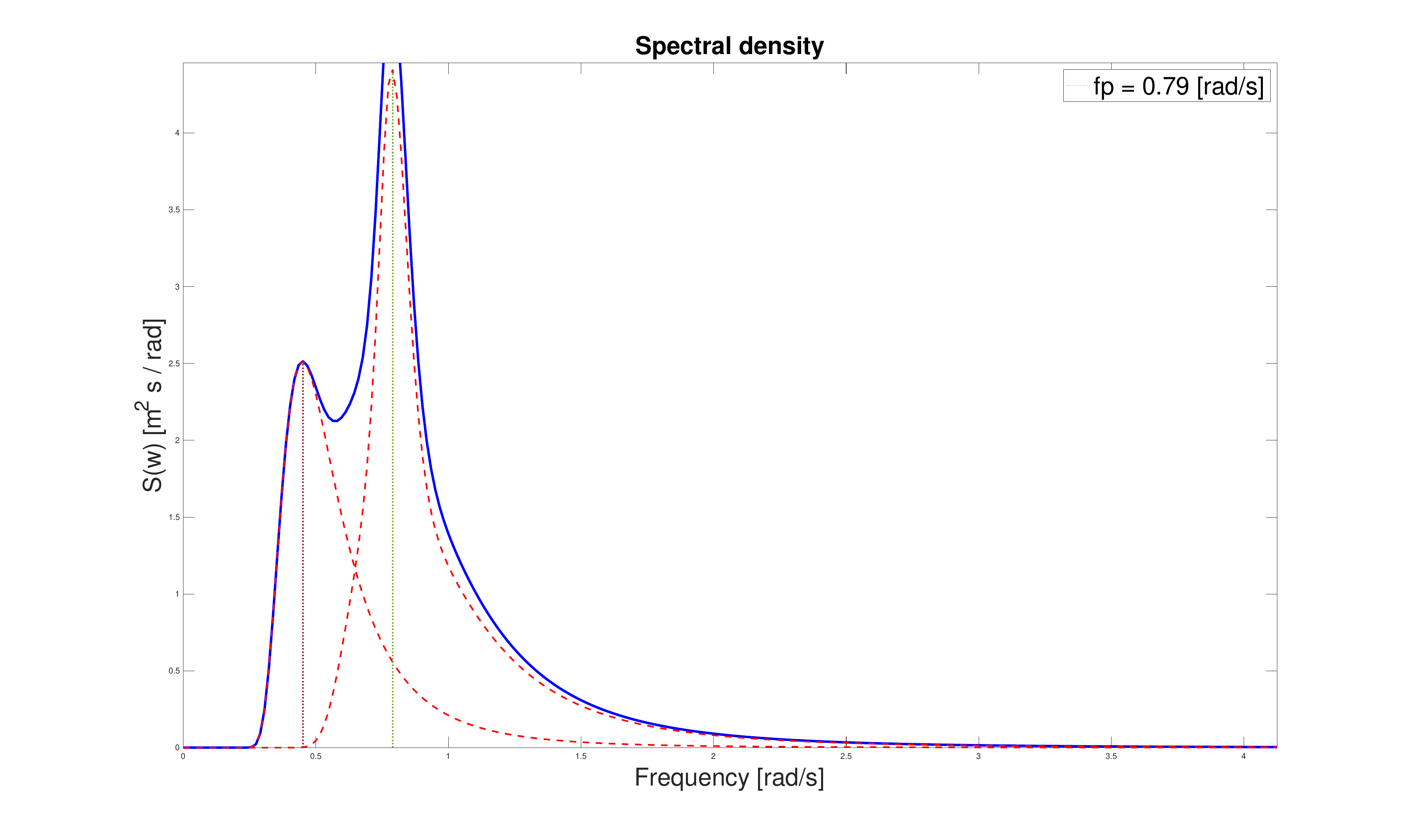}} \hspace{-10mm}
\makebox[0.48\textwidth][l]{\includegraphics[width=0.6\textwidth]{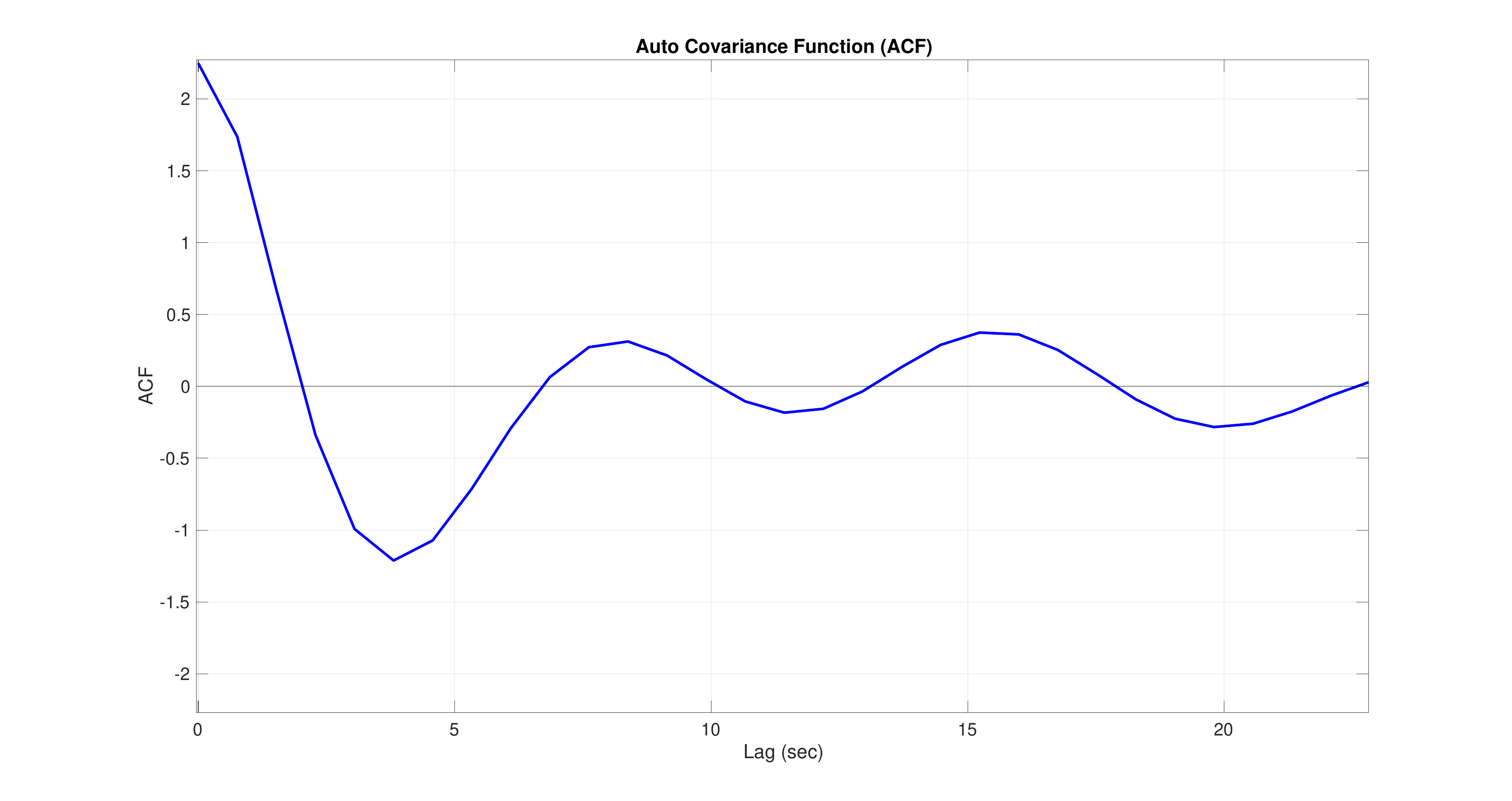}}\vspace{0cm}
\caption{\small The Torsethaugen spectrum {\it(left)} and the corresponding correlation {\it (right)} that have been used for generating data that illustrate the introduced concepts.}
\label{fig:spcov}
\end{figure}
The input to the model is significant wave height, $H_s$, for total sea and spectral peak period, $T_p$, for the primary (highest) peak. In our example, we have set $H_m = 6[m]$; $Tp  = 8[s]$; and used WAFO-toolbox \citep{gitwafo} to make all computations. In Figure~\ref{fig:spcov}, we see the graphs of the two-peaked spectrum and the corresponding autocovariance. Unless otherwise stated, a covariance is assumed to be normalized to an autocovariance throughout the rest of the paper. A simulated trajectory from the model is shown in blue in Figure~\ref{swcl} together with the clipped process in red.
In this example, the intervals will be dependent, which is generally the case for clipped processes, as mentioned previously. Therefore, the main idea behind the IIA is to overlook this complex dependency structure and approximate it with a binary process where the plus and minus interval distributions are independent and identically distributed (iid). Under this assumption, the relation between the plus and minus distribution directly and explicitly relates to the characteristics, such as the expected value or covariance. Therefore, the clipped process's characteristics can be used to derive an approximation of the excursion distribution through this simplified relation.

However, the key question is which characteristics of the clipped process should be used to construct this binary process. For example, the mean length of the intervals could be used and is given by the well-known formula from \cite{Rice}. However, since the expected value does not uniquely determine a distribution, a choice of distribution must still be made.

Another choice is to match the covariance of a stationary binary process to the one of the clipped. However, it follows from the inspection paradox that the interval containing zero of a stationary binary process will not have the same interval distribution as the other intervals. While the mathematical foundations and details have long been resolved, see \cite{Palm}, \cite{Khinchin}, and \cite{RyllNardzewski}, these slight differences need to be accounted for when matching the covariance to a stationary binary process.

A way to address this issue is to use a non-stationary binary process and match expected value functions instead. This is the main idea of this paper and the topic of Section~\ref{sec:slepianIIA}. Before presenting this approach, we introduce some important properties of the stationary and non-stationary binary processes used and some problems arising from using the covariance function in the ordinary IIA approach.

\vspace{2mm}
\raggedbottom
\subsection{The switch process} 
The switch process is a binary process similar to the alternating renewal process but takes the specific values one and minus one. This process is used to approximate the clipped process. Therefore, we present the two versions and some properties of the switch process in this section. We start with the non-stationary version.

\begin{figure}[t]
\rule{-1.8cm}{0cm}
{\includegraphics[width=0.53\textwidth]{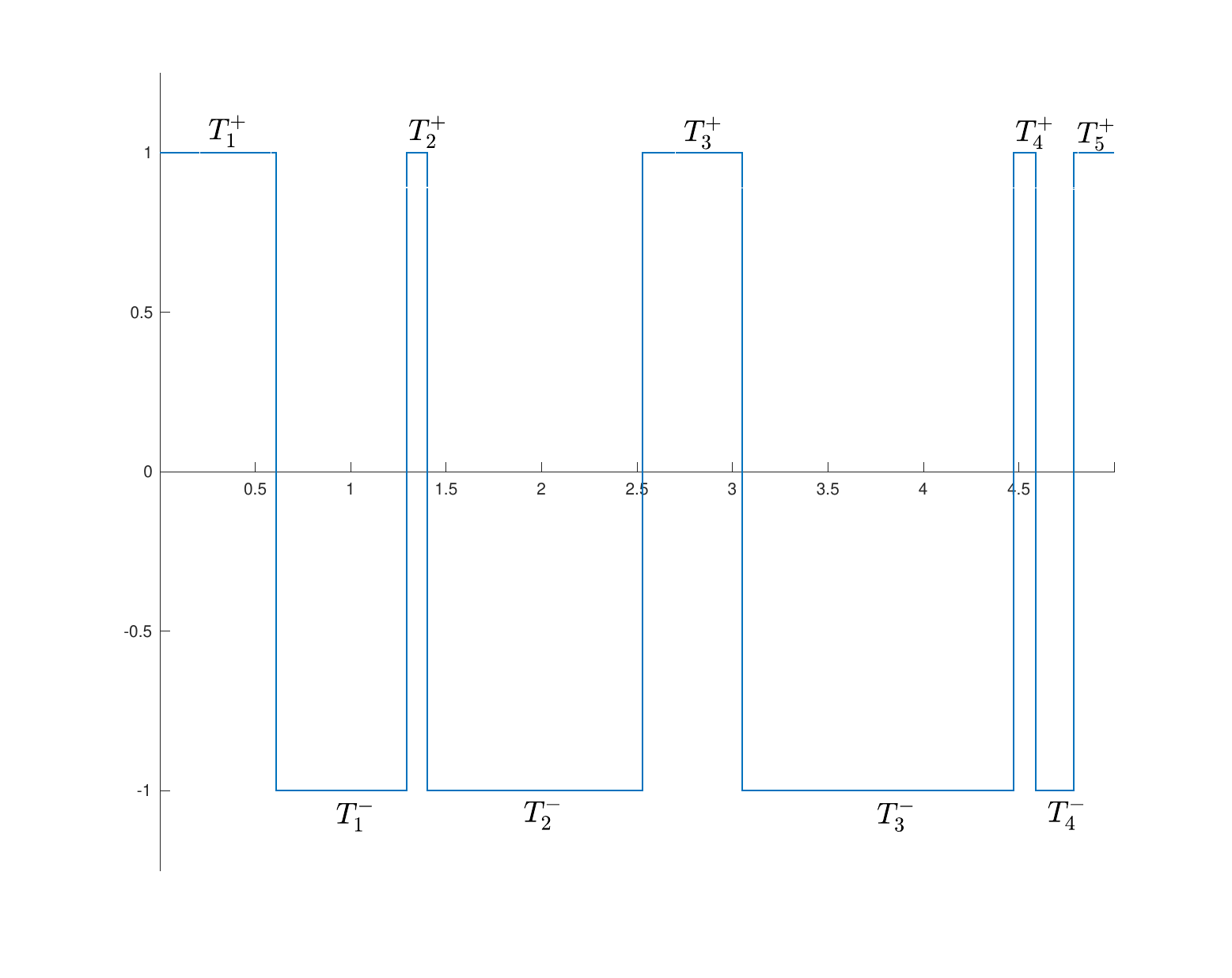}} \hspace{-10mm}
\makebox[0.4\textwidth][l]{\includegraphics[width=0.5\textwidth]{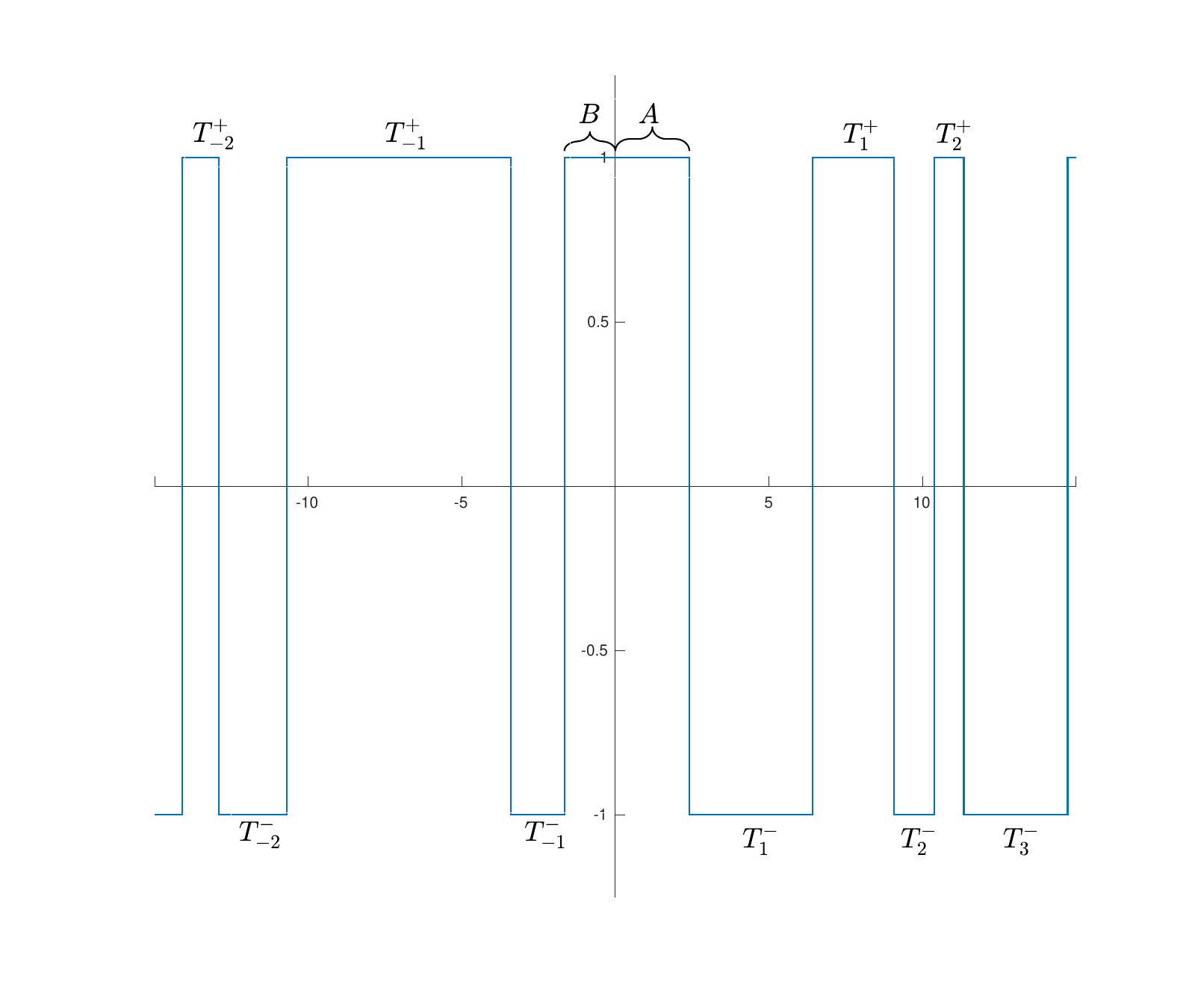}}\vspace{0cm}
\caption{\small The non-stationary switch process {\it(left)} and the stationary version, {\it (right)} with $\delta=1$. The inspection paradox states that the interval that contains zero has the distribution $A+B$, which is not the same as the other interval distributions.}
\label{fig:SP}
\end{figure}

Consider the positive line and its origin as a starting reference point and define a random process taking values one and minus one over interlaced intervals of the lengths $T_i^+$,  $T_i^-$, $i \in \mathbb N$. The distributions of these are referred to as the interval distributions to distinguish them from the excursion distribution. We put one over the interval $(0, T_1^+]$, after which we switch to minus one over the interval $(T_1^+, T_1^++T_1^-]$, then we switch again to one, and so on. More specifically,  
\begin{align}
    D(t)=(-1)^{N(t)},
\end{align} 
where $N(t)$ is the renewal process constructed from $T_i^+$ and $T_i^-$, $i\in \mathbb{N}$. The resulting one/minus one process is denoted by $D(t)$, $t>0$, and is referred to as the non-stationary switch process. In Figure \ref{fig:SP} (\textit{left}), a trajectory of such a process can be seen.

Since the excursion intervals, above and below zero, for a stationery Gaussian process are symmetric in distribution, we utilize this symmetry by letting $T_i^+\stackrel{d}{=}T_i^-$ for all $i\in\mathbb{N}$. 
Further, denote the common cdf of $T_i^+$ and $T_i^-$ by $F$ for which we assume that the density function $f$ exists and that the distribution has a finite expectation, i.e., ${\mu=\Ex T_i^\pm<\infty}$

Properties of $D(t)$, $t\ge 0$, can be determined from the distribution $F$ of $T_+$ and $T_-$. 
However, it is important to note that the process $D(t)$ is not stationary due to the special role the origin plays in its definition. For example, the expected value function $E(t)=\Ex{D(t)}$ is not constant with time.

Next, we present some key properties of the switch process, which was derived in \cite{Bengtsson2}. Since there is a wide use of the Laplace transform, denoted with $\mathcal{L}(\cdot)$, some properties of this transform are collected in Appendix \ref{App:Trans}. 
\begin{proposition} \label{prop:Esp}
Let $D(t)$ be a non-stationary switch process with the expected value function $E(t)=\Ex D(t)$. Then we have the relations for $s>0$ 
\begin{align}
\mathcal{L}(E)(s)&=\frac{1}{s}\frac{1-\Psi(s)}{1+\Psi(s)},
\\
\Psi(s)&=\frac{1-s\mathcal{L}(E)(s)}{1+s\mathcal{L}(E)(s)}.
\end{align}
where $\Psi(s)=\mathcal{L}(f)(s)$.
\end{proposition}
This proposition provides an explicit link between the expected value function and the distribution of the intervals. This connection will later form the basis for the Slepian-based IIA.

It's clear from the expected value function that the process $D(t)$ is not stationary. However, we can extend the process to the entire real line, use a standard tool from renewal theory, and delay it forward and backward around zero. Then, these delaying distributions can be chosen so that the resulting process is stationary. Since the value of this process at zero will be random, we also need a Bernoulli random variable to determine the value of the process on the interval that contains zero.

Let $A$, $B$ be the delays forward and backward and $\delta$ a Bernoulli random variable taking the value one and minus one. Then, we can define the delayed switch process by
\begin{align}\label{eq:SSP}
D_{s}(t)=
\begin{cases}
- \delta, & -B < t <A, \\
\delta D_+(t-A), &  t \geqslant A, \\
-\delta D_-(-(t+B)), & t \leqslant -B,  
\end{cases}
\end{align}
where $D_\pm$ are two independent non-stationary switch processes attached at the ends of the delaying interval $[-B, A]$ with the common interval distribution $F$. This construction is illustrated in Figure \ref{fig:SP} (\textit{left}). We now need the distribution of $A, B,\delta$ such that the delayed switch process becomes stationary. For this, we have the following proposition, which follows from applying the key renewal theorem.
\begin{proposition}
    \label{prop:stat}
    Let $D_+(t)$ and $D_-(t)$, $t>0$ be two independent non-stationary switch processes with the interval distribution $T$, with the distribution $F$ and density $f$, such that the expected value $\Ex T=\mu<\infty$. If $A,B,\delta$ are non-negative, mutually independent and independent of $D_\pm$ with a distribution is given by 
    \begin{align*}
        f_{A,B}(a,b)&=\frac{f(a+b)}{\mu} 
        \\
        f_A(t)&=f_B(t)=\frac{1-F(t)}{\mu}
        \\
        \Prob(\delta=1)&=\Prob(\delta=-1)=\frac{1}{2}
    \end{align*}
    Then, the delayed switch process in Equation (\ref{eq:SSP}) is stationary. We call a delayed switch process with this delay for the stationary switch process. 
\end{proposition}

\begin{remark}
    While it is not necessary that $T$ has a density, for the proposition to hold, it is only necessary that $T^+_1+T^-_1$ is not sitting on a lattice. 
\end{remark}

From Proposition \ref{prop:stat}, we observe the so-called inspection paradox. The interval containing zero will have a different distribution than the regular interval distribution $T$ since it will have the length $A+B$. We illustrate this in Figure \ref{fig:SP} (\textit{left}).

The main characteristic of the stationary switch process is the covariance function, which, for the symmetric case, when $T_+=T_-$ uniquely determines and characterizes the process. In the next proposition, we relate the covariance and the interval distribution.
\begin{proposition} \label{th:SwitchLapFour}
For a fixed $t\ge 0$, the distribution of $D_s(t)$ is uniquely characterized by $P_\delta(t)=\Prob(D_s(t)=1|D_s(0)=\delta)$, $\delta=\pm 1$ that have the following Laplace transforms
\begin{align*}
\mathcal L (P_\delta)(s)&=\frac{1}{s}
    \begin{cases} \displaystyle
        1-\frac 1 {\mu s} 
        \frac{ 1-\Psi(s)}{1+\Psi(s)} 
        &;\delta=1,
        \\
        \displaystyle
        \frac 1 {\mu s} 
        \frac{1-\Psi(s)}{1+\Psi(s)}&; \delta=-1,
    \end{cases}
\end{align*}
where $\Psi$ is the Laplace transform of the probability distribution of $T_i^{\pm}$'s. 
The covariance function $R(t)=\Cov{(D_s(h), D_s(t+h))}$ has the form
\begin{align*}
    R(t)=P_{1}(t)-P_{-1}(t)
\end{align*}
with the Laplace transform
\begin{align*}
    \mathcal L (R)(s)=\frac{2}{s\mu}\left(\frac{\mu}{2}- \frac1s \frac{1-\Psi(s)}{1+\Psi(s)} \right).
\end{align*}
\end{proposition}
This connection between the covariance and the interval distribution allows us to express the interval distribution as a function of the covariance in the following way
\begin{align} \label{eq:ordIIA}
    \Psi(s)=\frac{1-\frac{\mu s }{2} + \frac{\mu s^2}{2} \mathcal{L}(R)(s)}{1+\frac{\mu s }{2} - \frac{\mu s^2}{2} \mathcal{L}(R)(s)}
\end{align}
This equation agrees with formula (215) of \cite{BrayMS}. In the next section, we will see how this equation is used to obtain an approximation of the excursion distribution in the ordinary IIA framework. 

\vspace{2mm}
\raggedbottom
\subsection{The ordinary IIA}
In the previous section, we saw how the covariance characterizes the stationary process through its relation to the interval distribution. The ordinary IIA uses this relation to approximate this by substituting the covariances function in Equation (\ref{eq:ordIIA}) with the covariance of a clipped process; this is the essence of matching covariance.

An explicit relation exists between the covariance of a stationary Gaussian process and its clipped process for the zero level. The two are related by
\begin{align*}
    R_{cl}(t)=\frac{2}{\pi} \arcsin (r(t)). 
\end{align*}
and then from this, we directly obtain the ordinary IIA approximated excursion distribution by combining the above equation with Equation (\ref{eq:ordIIA}) in the following way 
\begin{align}
\label{eq:IIAcov}
    \Psi_{\text{IIA}}(s)=\frac{1-\frac{\mu s }{2} + \frac{\mu s^2}{2} \mathcal{L}(R_{cl})(s)}{1+\frac{\mu s }{2} - \frac{\mu s^2}{2} \mathcal{L}(R_{cl})(s)}. 
\end{align}
By doing this matching, we overlook the dependency structure and assume that it is sufficiently small to be ignored. However, two fundamental questions arise from this approach that have rarely been addressed in the past, see \cite{Lindgren_2022}.

The first question is: {\em  Does $\mathcal{L}(R_{cl})(s)$ in equation (\ref{eq:IIAcov}) always lead to a valid probability distribution, i.e., if the IIA is mathematically sound?} This question is equivalent to finding under what conditions a Laplace transform corresponds to probability density functions. There is Bernstein's Theorem, which states if $\Psi_{\text{IIA}}(s)$ is completely monotone, then it's the Laplace transform of a probability distribution, see \cite{Widder}, Theorem 12a, p.160. To show that a function is completely monotone involves all the derivatives of a function. This is difficult in practice. In fact, even determining when a rational function is a Laplace transform of a probability density has been so far an unsolved problem despite many partial results, see \cite{Zemanian, Zemanian61, SumitaM}

The second question, given an affirmative answer to the first one, is: {\em Is there any effective and explicit form of this distribution?} These questions will be approached, and partial answers will be obtained after an alternative approach is proposed. In this approach, one matches a non-stationary switch process with the clipped Slepian process with the origin attached to a zero-crossing, which is the topic of the next section.

\section{The Slepian-based IIA }
\label{sec:slepianIIA}
\noindent
The Slepian-based IIA approximates the excursion distribution by independent intervals when the reference point is attached to a level crossing instead of the origin. Conceptually, it is more natural because the crossing instant is the excursion-relevant event while the origin is not. It also avoids the origin location bias, also known as the inspection paradox, as discussed above through the variables $(A, B)$ in the section on the stationary switch process.

Compared to the ordinary IIA, the Slepian-based IIA uses different characteristics than the covariance. Instead of matching the covariances of the two stationary processes, it matches the expected value function of a switch process and the clipped Slepian process. 
The Slepian model was introduced by \cite{Slepian1963} and models the statistical behavior at the instants of $u$ level up-crossing of a stationary Gaussian process. 
This stationary Gaussian process needs to be sufficiently smooth to have well-defined level crossings and a covariance function that is at least twice continuously differentiable. For such a Gaussian process $X(t)$ with the covariance function $r(t)$ the Slepian process is given by 
\begin{align} \label{eq:Slepian}
    X_{u}(t)=u \cdot r(t) - R \cdot \frac{r'(t)}{\sqrt{-r''(0)}}+\Delta(t),
\end{align}
where $R$ is a standard Rayleigh variable, with the density $f_R(s)=se^{-s^2/2}$ independent of the non-stationary Gaussian process $\Delta$, with covariance 
\begin{align*}
    r_\Delta(t,s) = r(t-s)- r(t)r(s)+\frac{r'(t)r'(s)}{r''(0)}.
\end{align*}
Since we focus on the zero-level crossing, this process reduces to 
\begin{align} \label{eq:SlepianZero}
    X_{0}(t)=-R \cdot \frac{r'(t)}{\sqrt{-r''(0)}}+\Delta(t).
\end{align}
From the above equation, we see the importance of having finite crossing density since $R$ is scaled by $-r'(t)/\sqrt{-r''(0)}$. For more details on the Slepian process, see \cite{LLR83}.  With the Slepian process defined, we can now introduce the clipped Slepian process, which is defined by
\begin{align*}
    D_0(t)={\rm sgn}\left(X_0(t)\right), \, t\ge 0,
\end{align*}
where $X_0$ is defined in \eqref{eq:SlepianZero}. 
Since the Slepian process models an up-crossing, this equates to the clipped version starting from one. Hence, matching the clipped Slepian process with the switch process is a natural choice.

\subsection{The approximated excursion distribution}
Due to the lack of stationarity of both the clipped Slepian process and the switch process, it is more appropriate to match the expected value instead of covariance. The expected value function has an explicit form and is given in the following proposition. 
\begin{proposition} \label{Ecliped}
Let $E_0(t)=\Ex D_0(t)$ be the expectation of a Slepian process clipped at zero.
Then for $t\geqslant0$,
\begin{align} \label{eq:E0}
    E_0(t)=- \frac{1}{\sqrt{-r''(0)}} \frac{ r'(t)}{\sqrt{1-r(t)^2}},
\end{align}
where $r(t)$ is the covariance function of the original Gaussian process, for which it is assumed that $r''(t)$ is well defined for $t\geqslant0$. 
\end{proposition}
\begin{proof}
We have
\begin{align*}
    \Ex{(D_0(t))}=\Ex{\left(I(X_0(t)>0)-I(X_0(t)\leqslant0)\right)}=2\Prob(X_0(t)>0)-1,
\end{align*}
the probability of $X_0(t)$ being above zero, depends on the sign of $r'(t)$, and we have 
\begin{align*}
\Prob(X_0(t)>0)&=
    \begin{cases}
        \Prob\left(R<\Delta(t)\frac{\sqrt{-r''(0)}}{r'(t)}\right)&: r'(t) > 0, 
        \\
        \Prob(\Delta(t)>0) &: r'(t)= 0, 
        \\
        1- \Prob\left(R\le \Delta(t)\frac{\sqrt{-r''(0)}}{r'(t)}\right) &: r'(t) < 0.
    \end{cases} 
\end{align*}
The probability of a Rayleigh random variable being less than a normal random variable is easy to obtain, and regardless of the sign of $r'(t)$, we have
\begin{align*}
    \Prob(X_0(t)>0)=\frac{1}{2} \left( 1- \frac{1}{\sqrt{-r''(0)}} \frac{r'(t)}{\sqrt{1-r(t)^2}} \right),
\end{align*}
which leads to \eqref{eq:E0}.
\end{proof}
From the proposition, it's clear that the behavior of the expected value function $E_0$ of the clipped Slepian process is largely determined by the derivative of $r(t)$. This thus serves as a direct link to the covariance of the process we are interested in and the expected value function.

However, matching in the time domain is cumbersome, which is why it's done in the Laplace domain. This is done by substituting the expected value of the clipped Slepian process from Proposition \ref{Ecliped} into the second equation of Proposition \ref{prop:Esp}. Hence, we directly obtain an expression for the approximated excursion distribution in the Laplace domain. We summarize this in the next proposition.

\begin{proposition} \label{Prop:aproxdist}
    Let $E_0$ be the expected value of the Clipped Slepian process, then the Slepian-based IIA approximation of the excursion distribution denoted by $\widehat{T}$ is given in the Laplace domain by
    \begin{align*}
     \Psi_{\widehat{T}}(s)=\frac{1-s\mathcal{L}(E_0)(s)}{1+s\mathcal{L}(E_0)(s)}.
     \end{align*}
\end{proposition}
The approximated excursion distribution can be obtained in the time domain by inverting the above Laplace transform. This can be done either analytically or numerically. 
However, there are still several important questions that need to be answered. For example, what conditions on $E_0$ guarantee that the approximated excursion distribution is indeed a valid probability distribution? This is difficult to verify for arbitrary expected value functions due to Bernstein's theorem, as previously mentioned. Further, what is the relation between the Slepian-based IIA and the ordinary IIA for the approximating zero-level excursions? The next section deals with these questions and properties of the Slepian-based IIA.

\subsection{Properties of the Slepian-based IIA}
One of the main advantages of the Slepian-based IIA is that approximated excursion distribution can be retrieved explicitly for a large class of stationary Gaussian processes in terms of a stochastic representation. 
The following theorem on this representation thus serves as the probabilistic foundation for the Slepian-based IIA and follows directly from Theorem~1 and Corollary~2 in \cite{Bengtsson2}. 

\begin{theorem}\label{Th1}
If the expected value function $E_0(t)$ of the clipped Slepian process is differentiable with a non-positive derivative for $t\geqslant 0$. Then, the Slepian-based IIA approximation of the excursion distribution has the following stochastic representation
\begin{align}
\label{eq:geomsum}
   {\widehat{T}}=\sum_{i=1}^{\nu} \widetilde{T}_i.
\end{align}
Where $\widetilde{T}_i$ is an iid sequence of positive random variables having the survival function equal to $E_0(t)$ and $\nu$ follows a geometric distribution with the probability mass function $\Prob(\nu=k)=(1/2)^k$, $k\in \mathbb{N}$.
\end{theorem}
A random variable $\widehat{T}$, and its distribution with the stochastic representation of Theorem~\ref{Th1} is called $2$-geometric divisible, and the distribution of $\tilde T_i$'s is referred to as a geometric divisor of the distribution of $\widehat{T}$. For a brief introduction to this class of distributions, see \cite{Bengtsson2}. 
The representation of Theorem \ref{Th1} lacks usefulness unless there is a way to express $\widetilde{T}$ in terms of some characteristic. Luckily, there is such a way, which we highlight in the next corollary.
\begin{corollary} \label{rem:divisor}
    The cdf of the divisor $\widetilde{T}$ is given by 
    \begin{align*}
        F_{\widetilde{T}}(t)=1-E_0(t).
    \end{align*}
\end{corollary}
The implication of this is that there are several methods to obtain $\widehat{T}$ from $\widetilde{T}$. The choice of method can, therefore, be tailored depending on the particular properties of $\widehat{T}$ we are interested in. For example, the expected value of $\widehat{T}$ follows directly from Wald's equation and is $\Ex \widehat{T}=2\Ex \widetilde{T}$.

If the distribution of $\widehat{T}$ is sought, the standard method is through the Laplace transform since we deal with sums of random variables. In the Laplace domain, then we have the following relation between $\widehat{T}$ and $\widetilde{T}$
\begin{align}
    \label{eq:gd}
    \Psi_{\widehat{T}}=\frac{ \frac{1}{2}\Psi_{\widetilde{T}} }{ 1- \frac{1}{2} \Psi_{\widetilde{T} } }. 
\end{align}
If the analytical inverse Laplace transform of Equation (\ref{eq:gd}) is unwieldy, the representation still allows for numerical methods. It should also be noted that the stochastic representation in Theorem~\ref{Th1} allows for straightforward sampling of values from $\widehat{T}$, opening up the use of Monte Carlo methods to approach this problem.

We have seen how Theorem~\ref{Th1} provides the theoretical foundation for the Slepian-based IIA and has several practical consequences for obtaining the approximated excursion distribution. However, one question remains: what is the relation to the ordinary IIA framework?

The core idea of both the ordinary IIA and Slepian-based IIA are very similar, namely the idea of matching characteristics of a clipped process to a version of the switch process. 
It turns out that the two approaches coincide for the zero level, which we show in the following Theorem. 
\begin{theorem}
\label{equiv}
Let  $R(t)$ be the covariance of a stationary switch process and $E(t)$ be the expectation of the corresponding non-stationary switch process. Then 
\begin{align*}
    R'(t)=-\frac{2}{\mu}E(t),
\end{align*}
where $\mu$ is the first moment of the interval distribution.

Similarly, if $R_{cl}(t)$ is the covariance of a clipped smooth stationary Gaussian process and $E_0(t)$ is the expected value of the clipped Slepian process for the level $u=0$, then
\begin{align*}
    R_{cl}'(t)=-\frac{2}{\mu}E_0(t),
\end{align*}
where  $\mu=\frac{\pi}{\sqrt{-r''(0)}}$ represents also the average length of the excursion interval. 
\end{theorem}
\begin{proof}
Theorem~2 in \cite{Bengtsson2} establishes the first part of the result. For the second part 
we have 
\begin{align*}
    R_{cl}'(t)&=\frac{d}{dt}\left( \frac{2}{\pi}\arcsin{r(t)} \right)
    \\
    &
    =\frac{2}{\pi}\frac{r'(t)}{\sqrt{1-r(t)^2}}.
\end{align*}
Let $\mu=\frac{\pi}{\sqrt{-r''(0)}}$, then, by Proposition~\ref{Ecliped},
\begin{align*}
    R_{cl}'(t)&=-\frac{2}{\mu}E_0(t).
\end{align*}
Where $\mu$ is the average interval length of the excursion intervals is a well-known fact following from the ergodic theorem.
\end{proof}
\begin{remark}
    In Theorem \ref{equiv}, two seemingly unrelated binary processes have the same relation between expected value and covariance when the process is attached to a jump. It is an interesting question if this relation generalizes to all binary processes with finite jump intensity on any compact set such that the intervals between jumps are well defined. 
\end{remark}

We obtain the following rather unexpected equivalence between the two IIA approaches. 
\begin{corollary} 
Matching the covariance of the clipped process to the covariance of a stationary switch process is equivalent to matching the expected value of the clipped Slepian process with the expected value of a non-stationary switch process. 
\end{corollary}

Additionally, Theorem \ref{equiv} also has implications for Theorem \ref{Th1}. Since any covariance function $R(t)$ that leads to the clipped covariance function $R_{cl}(t)$ having a non-positive derivative for $t>0$ will correspond to an expected value function that satisfies the conditions of Theorem \ref{Th1}, the theorem also provides a foundation for the ordinary IIA framework for a large class of Gaussian processes.

We end this section with a caveat on when both the Slepian-based and ordinary IIA approaches encounter technical difficulties. 
If $\Ex \tilde T_i=\infty$ this implies that $\Ex \widehat{T}=\infty $ in (\ref{eq:geomsum}), which leads to a limitation for the choice of $E_0(t)$ in order to have a sensible IIA. 
\begin{corollary}
    If $S(t)$ is a survival function of a distribution such that ${\int_0^\infty S(t)~dt=\infty}$, then it cannot be the expected value function of any stationary switch process. 
    In particular, if a covariance function $r(t)$ of a Gaussian process yields a non-increasing and non-integrable expected value function 
    \begin{align*}
        E_0(t)=-\frac{1}{\sqrt{-r''(0)}}\frac{r'(t)}{\sqrt{1-r^2(t)}},
    \end{align*}
    then there is no switch process that can be used for the purpose of matching. 
\end{corollary}
With this caveat in place, we will next treat the general problem of estimating the persistency coefficients within the IIA framework.

\section{Persistency coefficient through the IIA} \label{sec:Montecarlo}
\noindent
The persistency coefficient is defined by \cite{BrayMS} as the asymptotic decay of the first persistence probability. This is the probability that a zero mean process $X(t)$, $t\geq0$ starting at zero does not change the sign between $(0,\tau)$. It has been argued in \cite{Sire2008} that if $\vert r(t)\vert < C/t$, then the tail of this probability decays exponentially.

In physics and engineering, the exponent, or coefficient, is an important characteristic of stochastic systems since it provides insight into how long a system can stay in a certain state. The problem is, as pointed out by \cite{BrayMS}, that this coefficient is process-specific and generally hard to find analytically. See \cite{BrayMS} for an overview of both the persistency coefficient and approximation methods for it. 

In the IIA framework, the persistency coefficient can be approximated by investigating the tail of $\widehat{T}$. If the conditions of Theorem \ref{Th1} are satisfied, this can be done through inverting Equation (\ref{eq:gd}) or the equation in Proposition \ref{Prop:aproxdist}. Since this is the form of a ratio of two functions, this can present some difficulty.

In physics, this difficulty has been circumvented by instead using the pole of Equation (\ref{eq:ordIIA}) with the largest negative real part as the approximation for the persistency coefficient.
There are several challenges to such an approach. Firstly, the Laplace transform considered needs to correspond to a probability distribution, which is not obvious, as it has already been discussed in previous sections.

Another challenge is whether finding a pole with the largest negative real value guarantees the corresponding exponential tail behavior. This can be justified if the distribution can be expressed in terms of a finite mixture of exponential distributions, \cite{SumitaM}, so its Laplace transform is a rational function. Alternatively, one can assume the rational form of the Laplace transform and determine the behavior by the partial fraction decomposition. However, for most stationary Gaussian processes used in physics, the covariance of the clipped process does not lead to a rational function. Moreover, the poles with the largest negative real part do not necessarily lay on the real axis.
In general, the mathematical justification of asymptotic behavior through the poles is not automatic. Thus, the use of the poles to find the persistency coefficient should be viewed as a heuristic one. However, when applicable, there is a simple criterion to find the pole within the Slepian-based IIA, which follows directly from Proposition \ref{Prop:aproxdist}. 
\begin{proposition}
\label{prop:pole0}
Let $\Psi_{\widehat{T}}$ be the Slepian-based IIA approximation from Proposition \ref{Prop:aproxdist}. Then, the largest negative real pole of $\Psi_{\widehat{T}}$ is found as the largest negative real solution to the equation
\begin{align*}
s\mathcal L (E_0)(s)+1=0.
\end{align*}
\end{proposition}

Additionally, if the conditions of Theorem \ref{Th1} are satisfied, then the largest negative real pole of the $\Psi_{\widehat{T}}$ can be found by investigating its divisor. The largest negative real pole of $\Psi_{\widehat{T}}$ corresponds to the largest negative real solution to the equation $\Psi_{\widetilde{T}}(s)=2$. It is clear that the tail of a geometrically divisible random variable is heavier than that of its divisor. However, there is no strict relation between the two tails, as the following example shows. 
Nevertheless, in the next proposition, we connect the exponential bounds of the tails.
\begin{example}
    Let $\tilde T_i = \alpha$, for $\alpha >0$, and $\nu$ geometric random variable with $p=1/2$. Then, by direct evaluation, $T=\sum_{i=1}^\nu \tilde T_i$ has the rate (persistency) $\log 2/\alpha$.
\end{example}
\begin{proposition}
Let $T$ have the stochastic representation of Theorem \ref{Th1} with the divisor $\widetilde{T}$ and assume that there is some $\beta>0$, then we have the two statements.   
\begin{itemize}
    \item[$i)$] If $\Prob(\tilde T >\tau) \leq e^{-\beta\tau},$ then $\Prob(T >\tau) \leq  e^{-\beta\tau/2}$.
    \item[$ii)$]  If $\Prob(\tilde T >\tau) \geq e^{-\beta\tau},$ then $\Prob(T >\tau) \geq  e^{-\beta\tau/2}$
\end{itemize}
\end{proposition}
\begin{proof}
The proof follows by induction, and we will prove $i)$ in full and $ii)$ can be obtained by changing the direction of the inequalities used in the proof. For $k\in \mathbb N$ independent copies $\tilde T_i$, we show that  
\begin{align*}
\Prob(\tilde T_1+\dots + \tilde T_k> \tau) \leq  e^{-b \tau}\sum_{j=0}^{k-1} \frac{b^j \tau^j}{j!}.
\end{align*}
Assume that the above is valid for a certain $k\in \mathbb N$, and $F_k$ is the cdf of $\tilde T_1+\dots + \tilde T_k$, then
\begin{align*}
\Prob(\tilde T_1+\dots + \tilde T_{k+1}> \tau)
&=
\int_0^\tau \Prob(\tilde T_{k+1} > \tau-s)~dF_k(s) + \Prob(\tilde T_1+\dots + \tilde T_k>\tau)
\\
&\leq 
\int_0^\tau e^{-b (\tau-s)}~dF_k(s) + \Prob(\tilde T_1+\dots + \tilde T_k>\tau)\\
&=
\left(\left(1-\Prob(\tilde T_1+\dots + \tilde T_k>s)\right)e^{-b (\tau-s)}\right) \Bigg|_{s=0}^\tau 
\\
&
-
\int_0^\tau \left(1-\Prob(\tilde T_1+\dots + \tilde T_k>s)\right)be^{-b (\tau-s)}~ds 
\\
&
+ \Prob(\tilde T_1+\dots + \tilde T_k>\tau)
\\
&= 1 - \int_0^\tau b e^{-b(\tau -s )}ds+
\int_0^\tau \Prob(\tilde T_1+\dots + \tilde T_k>s)be^{-b (\tau-s)}~ds\\
&=  e^{-b\tau}\left(1+b
\int_0^\tau \Prob(\tilde T_1+\dots + \tilde T_k>s)e^{b s}~ds\right)\\
&\leq 
e^{-b\tau}\left(1+\sum_{j=0}^{k-1} \frac{ b^{j+1}}{j!}
\int_0^\tau s^j~ds\right)
=e^{-b\tau}\left(1+\sum_{j=0}^{k-1} \frac{ b^{j+1}}{(j+1)!}
\tau^{j+1}\right),
\end{align*}
which concludes the induction proof. 
The result follows from 
\begin{align*}
\Prob(T\ge \tau)&=\sum_{k=1}^\infty \frac{1}{2^k} \Prob(\tilde T_1+\dots + \tilde T_k >\tau)
\leq 
   e^{-b \tau} \sum_{k=1}^\infty \frac{1}{2^k}\sum_{j=0}^{k-1} \frac{b^j \tau^j}{j!}
    =e^{-b \tau} \sum_{j=0}^\infty \frac{b^j \tau^j}{j!}\frac{1}{2^j}.
\end{align*}
\end{proof}

As stated previously, the relation between the tail of the divisor and the full distribution does not follow a general structure. This is why simulation-based methods might be preferable for approximating the persistency coefficient. 
Regardless if trajectories are simulated or samples are obtained using Theorem \ref{Th1}, there is still a need to approximate the persistency coefficient from these samples. This can be easily done under assumptions on the tail behavior.

Suppose that the tail has an exponential form, i.e., for large $t$, there is a $\theta>0$ such that 
\begin{align}
    \Prob({\widehat{T}}>t)\sim e^{-\theta t},
\end{align}
the persistency coefficient can be approximated using the ordinary least squares estimator for $\theta$ on the empirical survival function since 
\begin{align} \label{eq:log}
    \ln( \Prob({\widehat{T}}>t))\approx  a + \theta t.
\end{align}
Of course, more elaborate tail estimators could be used in this context, but since reasonably large samples are obtainable, striving for efficiency is not essential here. The core idea of estimating tail coefficients using linearization and least squares on the empirical distribution function or a QQ plot is well established. Early results of this idea were presented by \cite{expLS} (ECDF) and \cite{expQQ} for the QQ plot. However, in actual computation $\Prob({\widehat{T}}>t)$ is replaced by $1-\hat{F}(t)$, where $\hat{F}$ is the empirical CDF of the sample from $\widehat{T}$. This method alleviates the problems that follow from estimating persistency coefficients using the poles of the Laplace transform.

We conclude this section with a discussion of cases where the pole method is ill-suited to approximate the persistency coefficient and where the IIA simply fails to produce meaningful approximations. Not much is known about persistency coefficients when the covariance function oscillates. This is also one of the cases where using the pole method leads to problems. \cite{WilsonGauss} studied processes with the covariance function 
\begin{align*}
    r(t)=\cos(\alpha t) \ e^{-\frac{t^2}{2}},
\end{align*}
for which a large discrepancy was found between the largest negative real pole obtained using the IIA framework and the persistency coefficient obtained from simulations. It seems that care needs to be taken when applying the IIA framework to processes with oscillating covariance functions. This has direct consequences for practical application if one wants to model sea states using the Torsethaugen model, for example.

It is not only processes with oscillating covariance functions that pose a problem for the IIA. There is an entire class of processes where the IIA approach fails, which we illustrate with a counterexample. Recall the classical result of Theorem~2 from \cite{NewellR1962}.
\begin{theorem}
\label{th:tails}
If $X(t)$ is a Gaussian stationary process with $EX(t)= 0$ and $r(t)< Ct^{-\alpha}$ for some $\alpha> 0$, $C > 0$ and all $t > 0$, then for some $K>0$:
\begin{align*}
\Prob\left (X(t)\ge 0, t\in[0,\tau]\right) <
\begin{cases}
    e^{-K\tau}; & 1<\alpha,\\
    e^{-K\tau/\log \tau}; & \alpha=1,\\
    e^{-K\tau^\alpha}; &  0<\alpha<1.
\end{cases}
\end{align*}
\end{theorem}
From this theorem, we obtain the following lemma needed for the counterexample. 
\begin{lemma}
\label{th:counter}
 Let $X(t)$, $t\in \mathbb R$ be a stationary Gaussian process with twice differentiable covariance function, $r(t)$ such that for a certain $\beta>1$
\begin{align*}
 \lim_{t\rightarrow \infty} r'(t)=O(t^{-\beta}),
\end{align*}
and $X_0(t)$, $t\ge 0$ be the Slepian process associated with $X(t)$.
If the process satisfies also the assumptions of Theorem~\ref{Th1}, then 
\begin{align*}
\Prob(\widehat{T}\ge \tau)\ge O(\tau^{-\beta}),
\end{align*}
and there exists no $\beta'>1$ such that 
\begin{align*}
\Prob(X_0(t)>0, t\in (0,\tau))\ge O(\tau^{-\beta'}).
\end{align*}
\end{lemma}
\begin{proof}
We have 
\begin{align*}
 \Prob(\widehat{T}\ge \tau)&\ge \Prob(\tilde T \ge \tau)
 = - \frac{1}{\sqrt{r''(0)}}\frac{r'(\tau)}{\sqrt{1-r^2(\tau)}} 
=O\left(\tau^{-\beta} \right).
\end{align*}
For the second part, if there is $\beta'>1$ such that
\begin{align*}
\Prob(X_0(t)>0, t\in (0,\tau))\ge O(\tau^{-\beta'}),
\end{align*}
then 
\begin{align*}
    \Prob\left (X(t)\ge 0, t\in[0,\tau]\right) &=
    \frac{1}{2\mu} \int_\tau^\infty \Prob \left(X_0(t)>0, t\in (0,u)\right )~du
    \\
    &\ge O(\tau^{-(\beta'-1)},
\end{align*}
which contradicts Theorem~\ref{th:tails}. 
\end{proof}

From the above result, we see that for the processes satisfying Lemma~\ref{th:counter}, the tail behavior of $\widehat{T}$ based on the IIA is of the power order, and this cannot be the case for the actual distribution. Thus, under these circumstances, using the IIA to approximate the persistency coefficient or the tail behavior is not valid. The next example illustrates this.
\begin{example}
Consider a stationary Gaussian process with the covariance function 
\begin{align*}
    r(t)={\left(1+\frac{t^2}{2}\right)^{-\alpha}},\,\,\alpha>0.
\end{align*}
We note that this is a valid covariance function since it is the characteristic function of a generalized symmetric Laplace distribution. The corresponding covariance function of the clipped process is
\begin{align*}
    R_{cl}(t)=\frac{2}{\pi} \arcsin \left( \left(1+ \frac{t^2}{2}\right)^{-\alpha} \right). 
\end{align*}
and thus
\begin{align*}
    E_0(t)=\frac{ \sqrt{\alpha } t \left(1+{t^2}/{2}\right)^{-(1+\alpha)}}{\sqrt{1-\left(1+ {t^2}/{2}\right)^{-2\alpha}}}.
\end{align*}
This function satisfies the condition of Theorem \ref{Th1} and yields a divisor with finite expectation. Hence, the approximated excursion distribution $\widehat{T}$ will also have finite expectation. 
However, we observe that Theorem~\ref{th:tails} can be applied, and the tail behavior of $\widehat{T}$ has nothing to do with the tails of the excursion times.  
\end{example}

\section{Applications and examples}
\label{sec:examples}
\noindent
There are many covariance functions leading to positive and decreasing expected value functions of the clipped Slepian process. In these cases, Theorem~\ref{Th1} allows for both the approximation of the persistency coefficient, denoted by $\theta$ throughout this section, and the full approximate excursion distribution. A cautionary example is also presented to highlight the danger of naively using the IIA framework to approximate the persistency coefficient. We summarize the cases considered in this section in Table~\ref{Tab:RandE}. 
\begin{table}[h!]  
\caption{Processes and their characteristics.}
\label{Tab:RandE}
\footnotesize
\begin{tabular}{lll}\toprule
\sc Process $X(t)$ & \sc Covariance $r(t)$ & \sc Expected value $E_0(t)$\\ 
\midrule[0.9pt]
\multicolumn{3}{c}{}   \\
Diffusion in $\mathbb R^d$, $d\in \mathbb N$ & $\displaystyle \frac{1}{\cosh^{d/2}t/2}$ & $\displaystyle \frac{1}{\cosh t/2}\sqrt{\frac d2 \cdot \frac{\cosh^2 t/2-1}{\cosh^d t/2-1}}$\, \vspace{2mm}\\
Random acceleration  & $\left(3  -  e^{- \vert t \vert }\right)e^{-\vert t \vert/2 }/2 $ & $\displaystyle \sqrt{\frac{3}{4e^t-1}}$ \vspace{2mm}\\
Shifted Gaussian, $\alpha \in \mathbb R$ & $\cos \alpha t\cdot e^{-t^2/2}$ & $\displaystyle \frac{1}{\sqrt{1+\alpha^2}}  \frac{\alpha \sin \alpha t+t\cos \alpha t}{{\sqrt{1- \cos^2\alpha t \cdot e^{-t^2}}}}\cdot e^{-t^2/2}$  \vspace{2mm}\\
Matérn, $\nu \ge 2$  
& 
$ C_\nu \cdot t^\nu K_\nu \left(  {t} \right)$ 
& 
$\displaystyle \sqrt{2(\nu-1)}
\frac{t^{\nu} K_{\nu-1} \left( t \right)}{\sqrt{C_\nu^{-2}- {t} ^{2\nu} K_\nu^2 \left( {t} \right)}}
$, $C_\nu=\frac{2^{1-\nu}}{\Gamma(\nu)}$ 
\\
\bottomrule 
\end{tabular}
\end{table}
\subsection{The $d$-dimensional diffusion}
\begin{figure}[h]
    \centering
    \includegraphics[width=0.475\textwidth]{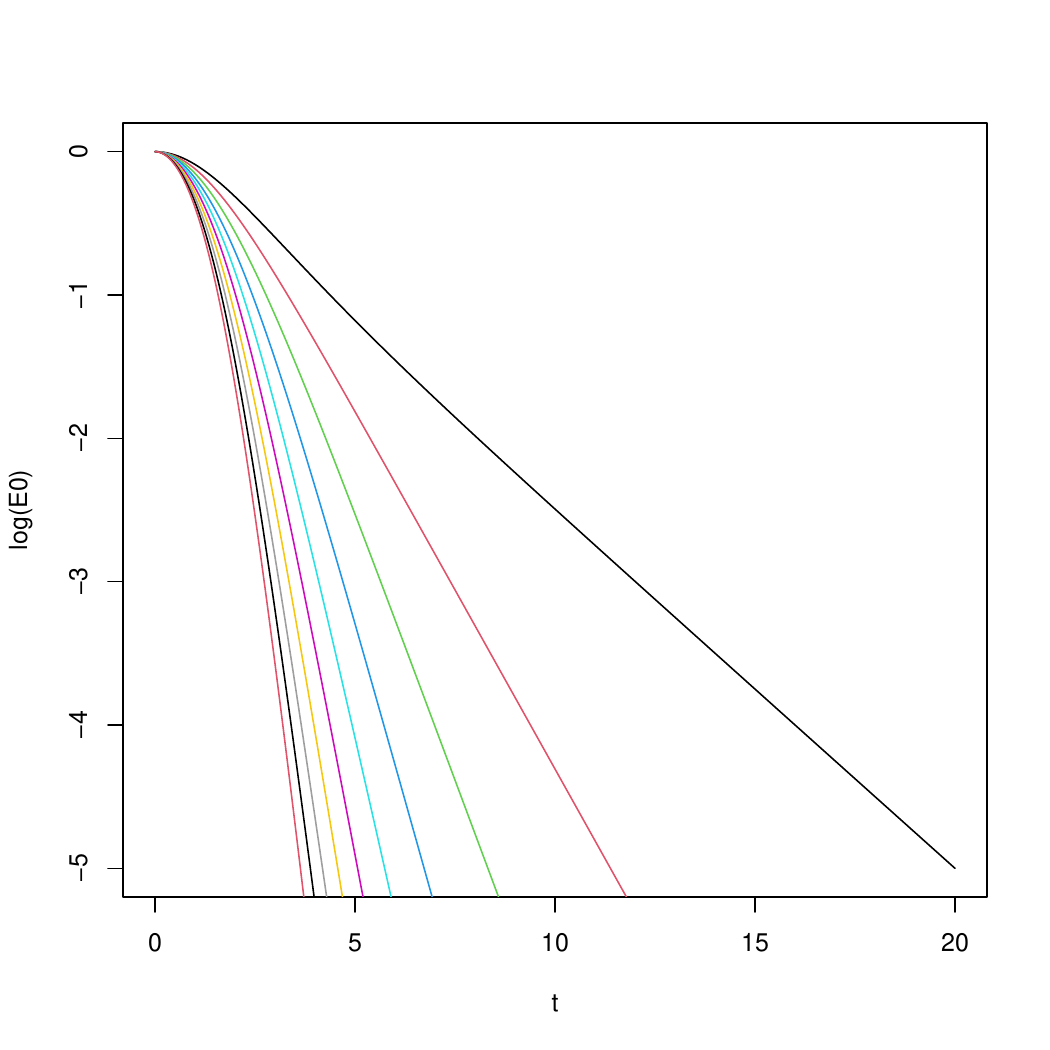}\hspace{-1mm}\raisebox{-3mm}{
    \includegraphics[width=0.48\textwidth]{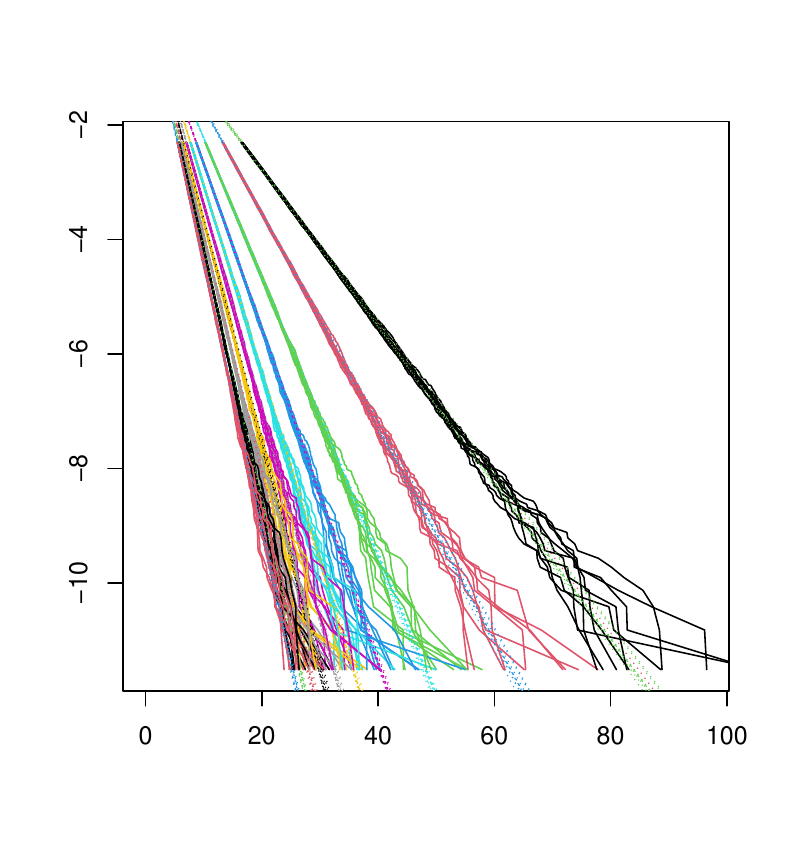}}\vspace{-9mm}
\caption{\small 
{\it (right)} The logarithm of $E_0$ for $d=1,...,10$. The steeper curves are associated with higher dimensions. 
{\it(left)} Logarithm of the empirical cdf from the approximated excursion distribution for the diffusion process, which is used for estimating the 
persistency coefficient.}
\label{fig:SP2}
\end{figure}

The $d$-dimensional diffusion is often referred to the scalar field $\phi (\mathbf x, t)$, with the location $\mathbf x \in \mathbb R^d$ and time $t>0$ that is the solution to 
\begin{align*}
\frac{\partial\phi}{\partial t}= \Delta \phi(\mathbf x,t)=\nabla^2\phi(\mathbf x,t),
\end{align*}
where the initial condition at $\tau=0$ is set to the $d$-dimensional (in argument) Gaussian noise.  
Then, for large $t$, the system is approximately described by 
a stationary Gaussian field
\begin{align*}
X(t)=\frac{\phi(\mathbf 0, e^t)}{\sqrt{\Ex (\phi^2(\mathbf 0, e^t)}},
\end{align*}
with the covariance function 
\begin{align*}
r(t)=\frac{1}{\cosh^{d/2}{(t/ 2)}}.
\end{align*}
The corresponding expected value function can be found in Table \ref{Tab:RandE} and it is clear that $E_0(t)$ is a strictly decreasing function in $t$ for each $d\in \mathbb N$ and thus satisfies the condition of Theorem~\ref{Th1}. It is also easy to observe that the tail of $E_0$ is almost exponential, which is illustrated in Figure \ref{fig:SP2} (\textit{left)}.

To obtain an approximation for the persistency coefficient, we simulate from the divisor using the outlined method in Appendix \ref{App:sim} and the stochastic representation of Theorem \ref{Th1}. The resulting persistency coefficient estimates can be found in Table \ref{tab:empstudy}. Each $95\%$ confidence interval is based on $10$ estimates of the persistency coefficient, which are in turn based on sampling $10^5$ observation from the divisor. The $ 10^4$ largest of these was then used for the least squares estimate of the persistency coefficients. In the table, we also include the results from a large and well-controlled simulation experiment of over $10^8$ realizations of first crossing events reported \cite{PhysRevLett.86.2712}. The exact computation method based on the generalized Rice formula and its implementation is given in the routing {\sc RIND} and presented in \cite{Lindgren_2022}. 

The conclusion from the simulation study is that the persistency coefficient obtained using the Slepian-based IIA is reasonable in comparison to the analytically derived values and the ones obtained by \cite{PhysRevLett.86.2712} and from using RIND. However, the method's limited accuracy should also be noted. 

\begin{table}[t!] 
\caption{Approximations of persistency coefficients for diffusions. \\For $d=1,2$, the true persistency is $0.1203$ and $0.1875$, respectively. }
\label{tab:empstudy}
\footnotesize
\centerline{
\begin{tabular}{rlll}\toprule
\sc Dimension &  Persistency coefficient & RIND2022 & NL2001 \\ 
\midrule[0.9pt]
1 & $0.1360\pm 0.0012$ & $0.1206 $ & $0.1205$ \\
2 & $0.1858\pm 0.0017$ & $0.1874$ & $0.1875$ \\
3 & $0.2441\pm 0.0014$ & $0.2382$ & $0.2382$ \\
4 & $0.2901\pm 0.0011$& $0.2805$ & $0.2806$ \\
5 & $0.3286\pm 0.0016$& $0.3171$ & $0.3173$ \\
6 & $0.3618\pm 0.0025$& * & * \\
7 & $0.3915\pm 0.0033$& * & * \\
8 & $0.4195\pm 0.0034$& * & * \\
9 & $0.4446\pm 0.0030$& * & * \\
10 & $0.4668\pm 0.0034$& $0.4589$ & $0.4587$ \\
\bottomrule 
\end{tabular}
}
\end{table}

\subsection{Random acceleration process}
The random acceleration process models a particle on the line whose acceleration is subject to a zero mean Gaussian noise force. In exponential time, i.e., $T=e^t$, the process is defined through the following equation of motion 
\begin{align*} 
    \frac{d^2Y}{dT^2}=\eta(T),
\end{align*}
where $\eta(T)$ is a white noise Gaussian process, i.e., the derivative of a Brownian motion (Wiener process) or a so-called $\delta$-correlated Gaussian white noise process, where $\delta$ stands for the Dirac's delta. More explicitly, $Y(T)=\int_0^T B(u)~du$, where $B$ is the standard Brownian motion. It is easy to check that the Gaussian process $X(t)=\sqrt{3/2}e^{-3t/2}Y(e^t)$ is a stationary Gaussian process with the covariance function 
\begin{align*}
    r(t)=\frac{3}{2}e^{-\frac{\vert t \vert}{2} }-\frac{1}{2}e^{-\frac{3\vert t \vert}{2} }.
\end{align*}
By elementary computation the first and second derivative of $r(t)$ are well defined with $r'(0)=0$ and $r''(0)=-3/4$. However, the third derivative has a discontinuity at zero, which \cite{Sire2007,Sire2008} argued leads to a poor approximation of the persistency coefficient. However, the process is sufficiently smooth since it has finite crossing intensity to construct a Slepian process. Hence, the expected value function of the clipped Slepian process for this example becomes 
\begin{align*}
    E_0(t)&=\frac{e^{-t/2}(1-e^{-t})\sqrt{3}}{\sqrt{4-(3-e^{-t})^2e^{-t}}}=\sqrt{\frac{3}{4e^t-1}}.
\end{align*}
This function satisfies the conditions in Theorem \ref{Th1}, and we use this observation to sample from the approximated excursion distribution through the divisor; see Appendix~\ref{App:sim}.

Taking a similar approach as for the diffusion example, we estimate the persistency coefficient based on a sample size of $10^7$ out of the $10^5$ largest where used for the least square estimate. This approach was repeated $10$ times which yielded the following $95\% $ confidence bound  
\begin{align*}
    \hat{\theta}=0.2647\pm 0.00083. 
\end{align*}
This approximation is in line with the estimates from \cite{Sire2007} ($\theta_{\text{app}}=0.2647$) using the pole method and which should be compared to the true value of $\theta=0.25$ derived by \cite{Sinai:1992aa}.
\vspace{2mm}
\raggedbottom
\subsection{Shifted and non-shifted Gaussian covariance}
Consider a stochastic process with the covariance function
\begin{align*}
    r(t)=\cos (\alpha t) \ e^{-{t^2}/{2}}, 
\end{align*}
for some $\alpha\geqslant0$. This case has been treated extensively by \cite{WilsonGauss}, which serves as a good overview of this case. We consider this example for three main reasons. The first is that it is the limiting case of the Mat\'ern covariance when the smoothness parameter goes to infinity. Secondly, the true persistency coefficient is not known for this process. Thirdly, for any $\alpha>0$, the covariance function is oscillating. Hence, we will first treat the case when $\alpha=0$ and then, after this, consider the case $\alpha=2$. 

The expected value function $E_0$ can be found in Table \ref{Tab:RandE} and for $\alpha=0$, it satisfies the conditions of Theorem \ref{Th1} and we use this observation when sampling from the approximated excursion distribution. The details can be found in Appendix \ref{App:sim}, and we sample $10^7$, out of which we use the $2 \cdot 10^6$ largest to estimate the persistency coefficient. This is repeated $10$ times, leading to the following $95\%$ confidence interval 
\begin{align*}
    \hat{\theta}=0.4116 \pm 0.00017. 
\end{align*}
To compare, we simulate trajectories using the WAFO package \citep{gitwafo}. The length of each trajectory is $10^7$, and $10^3$ is simulated. From these, we obtain an estimated persistency coefficient of $\theta_{\text{sim}}=0.4199\pm0.00058$. Comparing these two estimates to the IIA using the pole methods by \cite{Sire2008} ($\theta_{\text{IIA}}=0.4115$), we observe that they are all reasonably close.

Consider the case where $\alpha=2$. For this case, $E_0$ oscillates, and hence Theorem \ref{Th1}, can not be relied on to ensure the validity of the IIA approach. However, if one naively applies the IIA framework and uses the pole method to approximate the persistency coefficient, an estimate of $\theta_{\text{IIA}}=2.3522$ is obtained. Using the previously mentioned setup for simulating trajectories, an approximation of $\theta_{\text{sim}}=1.3795\pm 0.0012$ is obtained. Hence, it is clear that the use of the IIA and finding the pole is not suitable for approximating the persistency coefficient.

We only consider the largest negative real pole here. However, if a numerical search for the complex poles, a pair of conjugate poles is located at $1.05\pm  0.53 \ i$. While the real part of this pair is closer to the approximated persistency coefficient using simulation, it is still not an exact approximation. However, by multiplying this pair of complex conjugate poles, the value $1.3834$ is obtained, which is remarkably close to the value from simulated trajectories. This seems to hold for several $\alpha>0$. We, however, can't necessarily see the reason why these values are close to the persistency coefficients obtained by simulating trajectories. 
In conclusion, the case when $\alpha>0$ shows that care is needed since the IIA approach might break down for certain classes of covariance functions.

\subsection{Mat\'ern covariance}
\label{subsec:Matern}
\begin{figure}[t]
\includegraphics[width=0.33\textwidth]{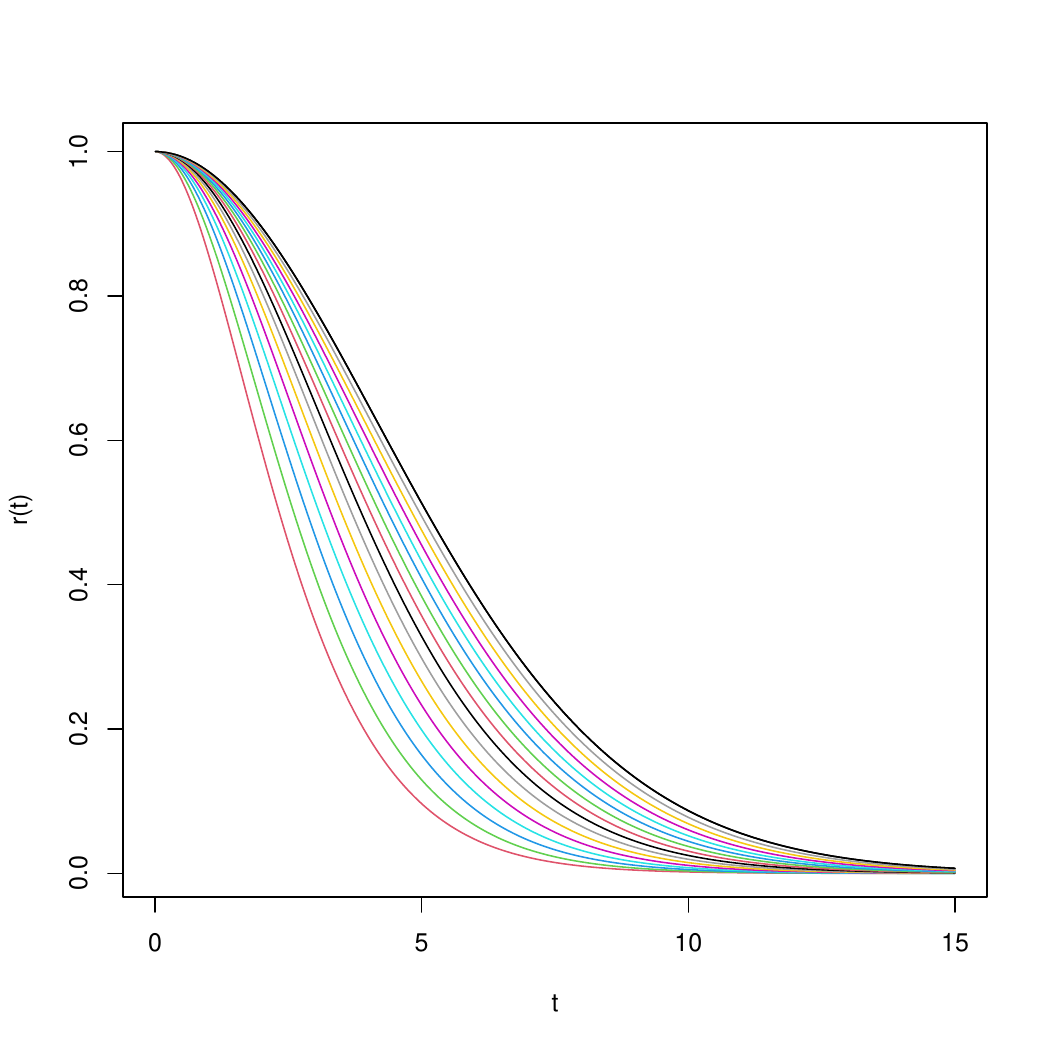}\includegraphics[width=0.33\textwidth]{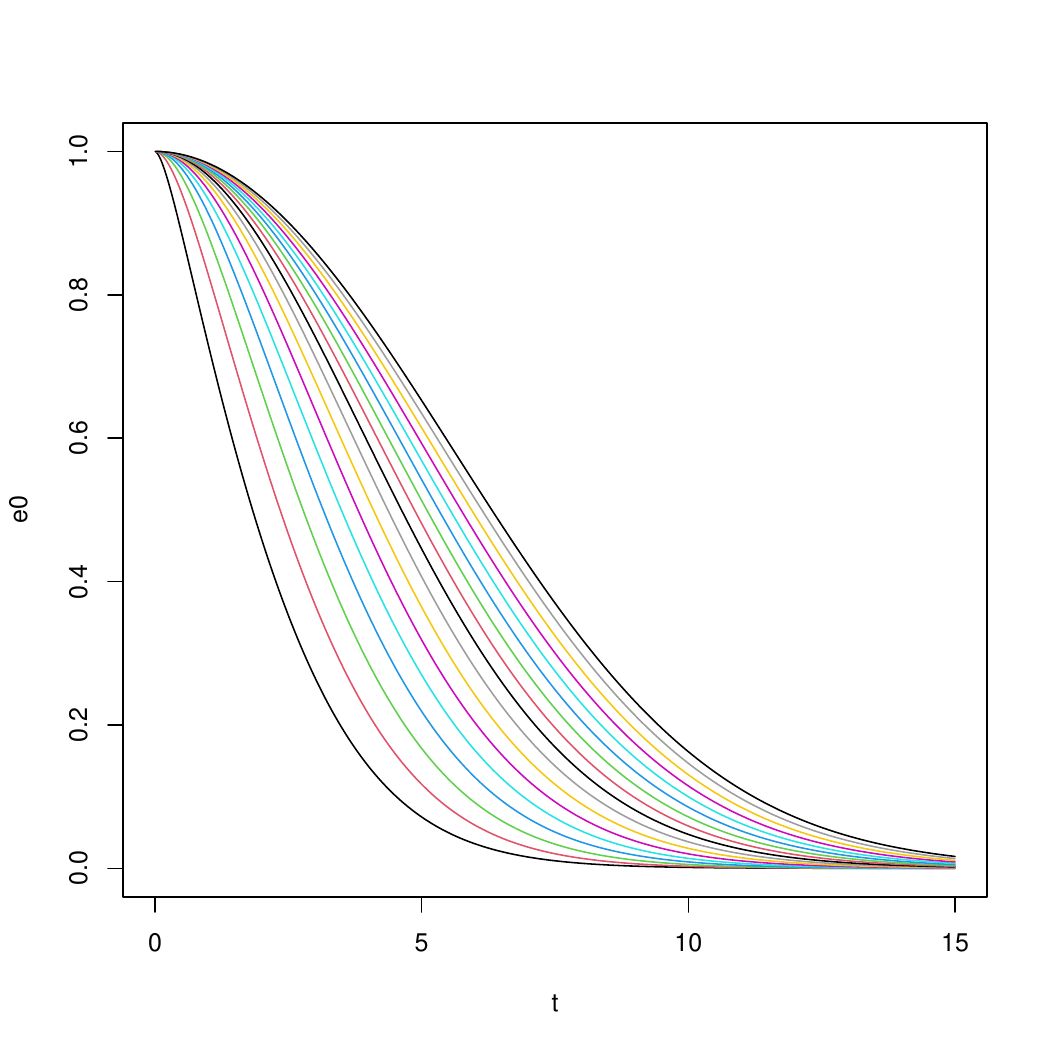}
\includegraphics[width=0.33\textwidth]{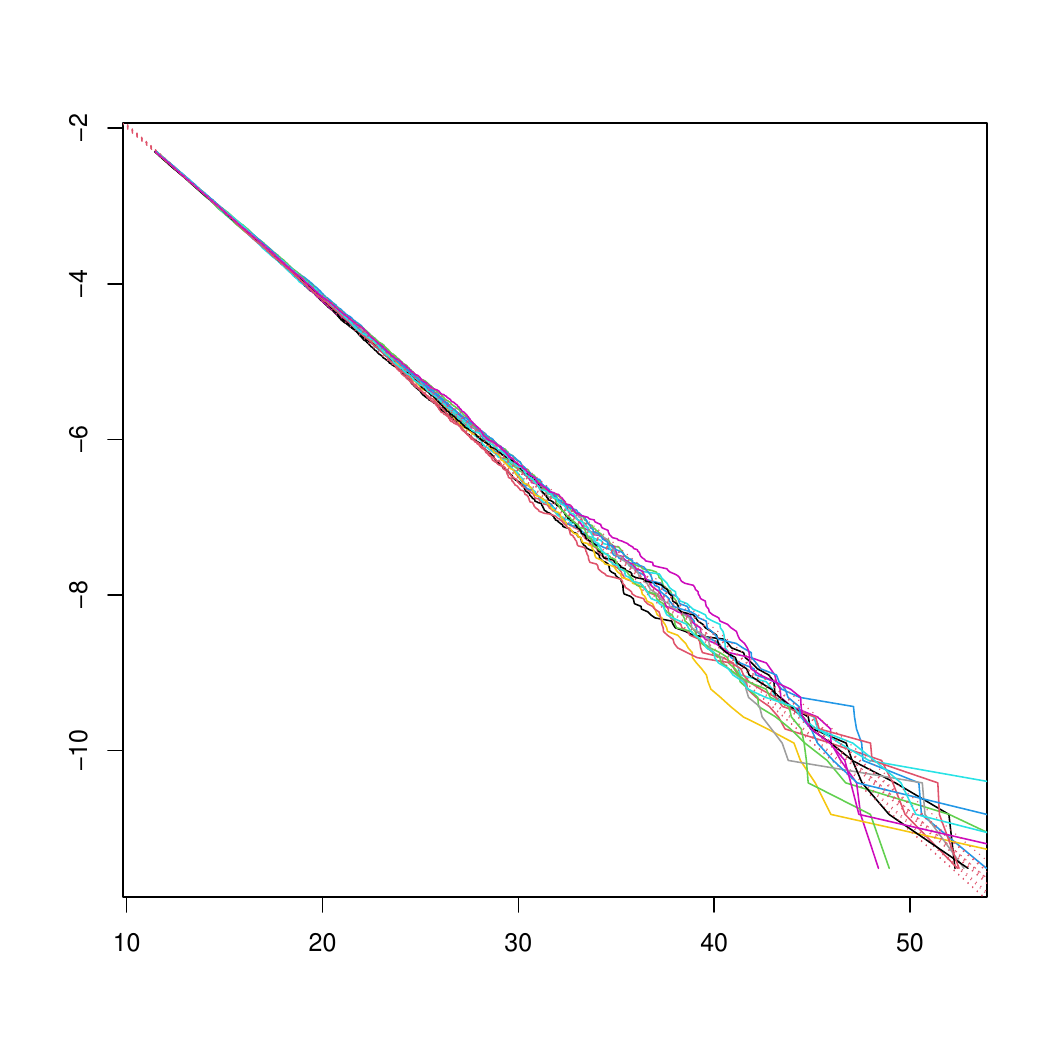}
\caption{The Mat\'ern covariance case. {\it Left:} The covariances $r$ for $\nu=l/2$, $l=4,5,\dots, 20$, the larger $\nu$ the stronger dependence. {\it Middle:} The survival functions $E_0$ of the geometric divisor, the larger $\nu$ the heavier tail.}
\label{Fig:Matern}
\end{figure}

The Matérn covariance 
\begin{align*}
r(t)=\frac{2^{1-\nu}}{\Gamma(\nu)} \cdot{t}^\nu K_\nu \left(  {t} \right),
\end{align*}
$\nu>0$ is widely used in spatial statistics. The parameter $\nu$ determines the smoothness of the process. Several well-known covariance functions are special cases of this class 
such as the exponential ($\nu=1$) and the previously treated Gaussian covariance ($\nu\rightarrow \infty$). In Figure~\ref{Fig:Matern} {\it (Left)}, the covariances of the underlying Gaussian processes are presented for selected values of $\nu$.

To ensure that the process has a finite crossing density, we limit ourselves to $\nu \ge 2$, thus allowing for the construction of the Slepian process. 
A key observation is that in this case, the expected value of the clipped Slepian process $E_0$ will satisfy the conditions of Theorem~\ref{Th1} and can be found in Table~\ref{Tab:RandE}. We present the survival functions of the divisor in Figure~\ref{Fig:Matern} {\it (Middle)} in which we see that the stronger dependence, as shown in the covariance $r$, translates to heavier tails in the divisor. 
\begin{remark}
     A re-parameterized version of the Matérn covariance is often used in spatial statistics. By letting $t=\sqrt{2\nu} \ d/\rho$, where $\rho$ is an additional parameter and $d\geqslant0$ is the new argument. The link between the persistency coefficients in the two parameterizations can be derived by observing that if $P(T>t)\approx C \ e^{-\theta t}$, then $P(T>\sqrt{2\nu} \ d/\rho) \approx C \ e^{-\theta \sqrt{2\nu} \ d/\rho}$. The re-parameterized persistency coefficient can be obtained by simply scaling the estimate by $\sqrt{2\nu} \ d/\rho$. 
\end{remark}
Generating random variables from these distributions can be numerically implemented, although the presence of the modified Bessel function may constitute a numerical challenge. For the special cases, $\nu = k+1/2$, $k\in\mathbb N$, we have the explicit formulas for the Bessel functions, and this can be utilized in random number generations, as discussed in Appendix \ref{App:sim}.

In Figure~\ref{Fig:Matern}~{\it (Right)}, the results of simulations in the case of $\nu=5/2$ are presented. This value was chosen because the random simulation from the geometric divisor could be devised by a simple inversion of $E_0$. Simulations for this case lead to the IIA persistency coefficient $ \hat{\theta}=0.2188\pm 0.0011$. The methodology of the obtained approximation has been described in previous sections. This is remarkably close to the persistency coefficient derived using WAFO. Using the same method for WAFO simulations as before, we obtain a $95\%$ interval of $\theta_{\text{sim}}=0.2184\pm 0.00026$.

\section{Conclusion}
\noindent
The IIA method of obtaining the excursion distributions has been explored, and its benefits and limitations have been determined. Several important questions have been posed, such as when the IIA framework is mathematically sound and whether explicit expressions of the approximated excursion distribution exist. By introducing the clipped Slepian process and matching its expected value function with that of the switch process, the validity of the IIA framework is shown for a large class of covariance functions. For this class, an explicit representation of the approximated excursion distribution exists. This representation allows additional ways to estimate persistency coefficients using Mote Carlo methods. A natural extension of this approach to non-zero level crossings will be investigated in the future. 

\vspace{-5mm}
\section{Acknowledgement}
\noindent
Henrik Bengtsson acknowledges financial support from the Royal Physiographic Society in Lund. Krzysztof Podg\'orski acknowledges the financial support of the Swedish Research Council (VR) Grant DNR: 2020-05168.

\bibliographystyle{apalike}
\bibliography{References}

\begin{thebibliography}{}

\bibitem[Bachelier, 1900]{Bachelier}
Bachelier, L. (1900).
\newblock Th\'eorie de la sp\'eculation.
\newblock {\em Annales scientifiques de l'\'Ecole Normale Sup\'erieure}, 3e s{\'e}rie, 17:21--86.

\bibitem[Bengtsson, 2024]{Bengtsson2}
Bengtsson, H. (2024).
\newblock Characteristics of the switch process and geometric divisibility.
\newblock {\em Journal of Applied Probability}, 61(3):802–809.

\bibitem[Brainina, 2013]{Brainina2013}
Brainina, I. (2013).
\newblock {\em Applications of random process excursion analysis}.
\newblock Elsevier, London.

\bibitem[Bray et~al., 2013]{BrayMS}
Bray, A., Majumdar, S., and Schehr, G. (2013).
\newblock Persistence and first-passage properties in nonequilibrium systems.
\newblock {\em Advances in Physics}, 62(3):225--361.

\bibitem[Khinchin, 1955]{Khinchin}
Khinchin, A. (1955).
\newblock Mathematical methods of the theory of mass service.
\newblock {\em Trudy Mat. Inst. Steklov.}, 49.

\bibitem[Kratz and Resnick, 1996]{expQQ}
Kratz, M. and Resnick, S. (1996).
\newblock The qq-estimator and heavy tails.
\newblock {\em Communications in Statistics. Stochastic Models}, 12(4):699--724.

\bibitem[Leadbetter et~al., 1983]{LLR83}
Leadbetter, M.~R., Lindgren, G., and Rootzen, H. (1983).
\newblock {\em Extremes and Related Properties of Random Sequences and Processes}.
\newblock Springer Series in Statistics Ser. Springer New York, New York, NY, 1st ed edition.

\bibitem[Lindgren et~al., 2022]{Lindgren_2022}
Lindgren, G., Podg{\'o}rski, K., and Rychlik, I. (2022).
\newblock Effective persistency evaluation via exact excursion distributions for random processes and fields.
\newblock {\em Journal of Physics Communications}, 6(3):035007.

\bibitem[Longuet-Higgins, 1962]{LonguetHiggins1962}
Longuet-Higgins, M. (1962).
\newblock The distribution of intervals between zeros of a stationary random function.
\newblock {\em Philos. T. Roy. Soc. A}, 254:557--599.

\bibitem[Majumdar et~al., 1996]{MajumdarSBC}
Majumdar, S.~N., Sire, C., Bray, A.~J., and Cornell, S.~J. (1996).
\newblock Nontrivial exponent for simple diffusion.
\newblock {\em Phys. Rev. Lett.}, 77:2867--2870.

\bibitem[McFadden, 1956]{McFadden1956}
McFadden, J. (1956).
\newblock The axis-crosssing intervals of random functions.
\newblock {\em IRE Trans. Inf. Theory}, IT-2:146--150.

\bibitem[McFadden, 1957]{McFadden1957AMS}
McFadden, J. (1957).
\newblock The variance of zero-crossing intervals [abstract].
\newblock {\em Ann. Math. Statist.}, 28:529.

\bibitem[McFadden, 1958]{McFadden1958}
McFadden, J. (1958).
\newblock The axis-crosssing intervals of random functions -- {II}.
\newblock {\em IRE Trans. Inf. Theory}, IT-4:14--24.

\bibitem[Molchan, 1999]{Molchan}
Molchan, G.~M. (1999).
\newblock Maximum of a fractional brownian motion: Probabilities of small values.
\newblock {\em Communications in Mathematical Physics}, 205(1):97--111.

\bibitem[Newell and Rosenblatt, 1962]{NewellR1962}
Newell, G.~F. and Rosenblatt, M. (1962).
\newblock Zero crossing probabilities for gaussian stationary processes.
\newblock {\em The Annals of Mathematical Statistics}, 33(4):1306--1313.

\bibitem[Newman and Loinaz, 2001]{PhysRevLett.86.2712}
Newman, T.~J. and Loinaz, W. (2001).
\newblock Critical dimensions of the diffusion equation.
\newblock {\em Phys. Rev. Lett.}, 86:2712--2715.

\bibitem[Palm, 1943]{Palm}
Palm, C. (1943).
\newblock Intensitatsschwankungen im fernsprechverkehr.
\newblock {\em Ericsson Technics}, 44:1--89.

\bibitem[Poplavskyi and Schehr, 2018]{PoplavskyiSchehr2018}
Poplavskyi, M. and Schehr, G. (2018).
\newblock Exact persistence exponent for the 2d-diffusion equation and related {K}ac polynomial.
\newblock {\em Phys. Rev. Lett.}, 121:1506011--1506017.

\bibitem[Rainal, 1962]{Rainal1962}
Rainal, A. (1962).
\newblock Zero-crossing intervals of {G}aussian processes.
\newblock {\em IRE Trans. Inf. Theory}, IT-8:372--378.

\bibitem[Rainal, 1963]{Rainal1963}
Rainal, A. (1963).
\newblock Zero-crossing intervals of random processes.
\newblock Technical Report AF-102, The Johns Hopkins University, Carlyle Barton Laboratory, Baltimore, MD.

\bibitem[Rice, 1944]{Rice}
Rice, S.~O. (1944).
\newblock The mathematical analysis of random noise.
\newblock {\em Bell Syst Tech J}, 23(3):282--332.

\bibitem[Ryll-Nardzewski, 1961]{RyllNardzewski}
Ryll-Nardzewski, C. (1961).
\newblock Remarks on processes of calls.
\newblock In of~Calif.~Press, U., editor, {\em Proc. Fourth Berkeley Symp. on Math. Statist. and Prob.}, volume~2, pages 455--465.

\bibitem[Schehr and Majumdar, 2007]{PhysRevLett.99.060603}
Schehr, G. and Majumdar, S.~N. (2007).
\newblock Statistics of the number of zero crossings: From random polynomials to the diffusion equation.
\newblock {\em Phys. Rev. Lett.}, 99:060603.

\bibitem[Schultze and Steinebach, 1996]{expLS}
Schultze, J. and Steinebach, J. (1996).
\newblock On least squares estimates of an exponential tail coefficient.
\newblock {\em Statistics \& Risk Modeling}, 14(4):353--372.

\bibitem[Sinai, 1992]{Sinai:1992aa}
Sinai, Y.~G. (1992).
\newblock Distribution of some functionals of the integral of a random walk.
\newblock {\em Theoretical and Mathematical Physics}, 90(3):219--241.

\bibitem[Sire, 2007]{Sire2007}
Sire, C. (2007).
\newblock Probability distribution of the maximum of a smooth temporal signal.
\newblock {\em Phys. Rev. Lett.}, 98:020601.

\bibitem[Sire, 2008]{Sire2008}
Sire, C. (2008).
\newblock Crossing intervals of non-{M}arkovian {G}aussian processes.
\newblock {\em Phys. Rev. E}, 78:011121--1--21.

\bibitem[Slepian, 1963]{Slepian1963}
Slepian, D. (1963).
\newblock On the zeros of gaussian noise.
\newblock {\em Time series analysis}, page 104–115.

\bibitem[Sumita and Masuda, 1987]{SumitaM}
Sumita, U. and Masuda, Y. (1987).
\newblock Classes of probability density functions having laplace transforms with negative zeros and poles.
\newblock {\em Advances in Applied Probability}, 19(3):632--651.

\bibitem[Torsethaugen and Haver, 2004]{Torsethaugen}
Torsethaugen, K. and Haver, S. (2004).
\newblock Simplified double peak spectral model for ocean waves.
\newblock In {\em {Proceedings of the 14th International Offshore and Polar Engineering Conference}}.

\bibitem[WAFO-group, 2017]{gitwafo}
WAFO-group (2017).
\newblock Wafo-project.
\newblock \url{https://github.com/wafo-project}.

\bibitem[Widder, 1946]{Widder}
Widder, D. (1946).
\newblock {\em The Laplace transform}.
\newblock Princeton University Press, Princeton.

\bibitem[Wilson and Hopcraft, 2017]{WilsonGauss}
Wilson, L. R.~M. and Hopcraft, K.~I. (2017).
\newblock Periodicity in the autocorrelation function as a mechanism for regularly occurring zero crossings or extreme values of a {Gaussian} process.
\newblock {\em Physical Review E}, 96(6):062129.

\bibitem[Zemanian, 1959]{Zemanian}
Zemanian, A.~H. (1959).
\newblock On the pole and zero locations of rational { Laplace } transformation of non-negative functions.
\newblock {\em Proceedings of the American Mathematical Society}, 10:868--872.

\bibitem[Zemanian, 1961]{Zemanian61}
Zemanian, A.~H. (1961).
\newblock On the pole and zero locations of rational { Laplace } transformation of non-negative functions {I}{I}.
\newblock {\em Proceedings of the American Mathematical Society}, 12(6):870--874.

\end{thebibliography}

\newpage
\appendix
\section{Laplace transform}
\label{App:Trans}
\noindent
Here, we establish the notation of the Laplace, which is the main tool used in many of our arguments.
Let us recall that for a positive and bounded function $P(t)$, $t>0$, its Laplace transform $\mathcal L (P)(s)$, $s>0$, is defined through 
\begin{align*}
\mathcal L (P)(s)=\int_0^\infty P(t)e^{-ts} ~dt.
\end{align*}
The Laplace transform can also be extended for $s<0$, and since it is a decreasing function of $s$, then there exists the smallest number $s_0\ge -\infty$ for which $\mathcal L (P)(s)$ is finite for all $s>s_0$. This value is referred to as the largest negative pole of the respective Laplace transform. For a distribution on the positive half line given by a cdf $F$, its Laplace transform $\Psi$ is given by
\begin{align*}
\Psi(s)=\int_0^\infty e^{-ts} ~dF(t).
\end{align*}
so that if the distribution has a density $f$, then
\begin{align*}
\mathcal L (f)(s)=\Psi(s).
\end{align*}
Similarly, we have the Laplace transform of the CDF
\begin{align*}
    \label{cdf}
\mathcal L (F)(s)=\Psi(s)/s.
\end{align*}

\section{Sampling from the Geometric divisor}
\label{App:sim}
\noindent
Here, we provide details of the simulators of various geometric divisors utilized in the IIA approximations based on Theorem~\ref{Th1}.

\subsection{Diffusion}
\noindent
Let us start from the diffusion in $d=2$. 
Since ${\rm arccosh}(x)= \\ { \log(x+\sqrt{x^2-1})}$, $x>1$, we have the following generator of the geometric divisor
\begin{equation}
    \label{eq:divdif2}
\tilde T_2=2\ln\left(1+\sqrt{1-U^2} \right) - 2 \ln U.
\end{equation}
For the diffusion in dimension $d=1$, the survival function can be written as 
\begin{align*}
E_{0,1}(t)= \frac{\sqrt{2}}2\sqrt{\cosh^{-1} (t/2) + 
 \cosh^{-2}(t/2)},
\end{align*}
which has the following survival inverse simulator
\begin{align*}
\tilde T_1=2\ln\left(1+\sqrt{1+8U^2}+\sqrt{2}\sqrt{1+4U^2-\sqrt{1+8U^2}}\right)-4(\ln U +\ln 2).
\end{align*}
For $d>2$, we have 
\begin{align*}
E_{0,d}(t)=E_{0,d-1}(t)\sqrt{\frac{d(\cosh^{d-1}(t/2)-1)}{(d-1)(\cosh^d(t/2)-1)}}.
\end{align*}
It is easy to verify that 
\begin{align*}
G_d(t)=\frac{d(\cosh^{d-1}(t/2)-1)}{(d-1)(\cosh^d(t/2)-1)}
\end{align*}
is also a survival function; thus if one can obtain its inverse $G_d^{-1}$, then the recursive formula for a sampler from $E_{0,d}$ is 
\begin{align*}
\tilde T_d=\min\left( \tilde T_{d-1},G_d^{-1}(U^2)\right),
\end{align*}
where the uniform variable $U$ is independent of $\tilde T_{d-1}$. To find the inverse $G_d^{-1}$, let us consider the function $b_d(a)$, $a>1$, $d\ge 3$ that is the inverse to the convex polynomial function on the positive half-line
\begin{align*}
a_d(b)=1+b+\dots + b^{d-1}.
\end{align*}
Then 
\begin{align*}
G_d^{-1}(g)=2\ln\left(1/{b_d\left(\frac{d}{d-g(d-1)}\right)}+\sqrt{1/b_d^2\left(\frac{d}{d-g(d-1)}\right)-1}\right).
\end{align*}
Since $a_d(b)$, $b>0$ is a convex polynomial function, one can efficiently evaluate its inverse using, for example, the Newton–Raphson method.

\subsection{Random acceleration process}
\noindent
Obtaining a random number from the divisor for the random acceleration process is straightforward by inverting the survival function given by $E_0$, leading to
\begin{align*}
\tilde T=\ln(2/U^2+1)-2\ln 2.
\end{align*}
This allows for fast simulation of the divisor and, hence, the full approximated excursion distribution through the stochastic representation of Theorem \ref{Th1}.

\subsection{Shifted Gaussian}
\noindent
We consider only the Gaussian case $(\alpha=0)$ in this section since Theorem \ref{Th1} does not apply to the shifted Gaussian covariance function. Compared to the two previous examples, the Gaussian case does not allow for easy inversion sampling for the divisor. However, the target density is explicit and can be evaluated at a low computational cost. Hence, random numbers can be simulated using rejection sampling. 
The target density denoted by $f$ is 
\begin{align*}
    -\frac{d}{dt}E_0(t)=\frac{e^{t^2}(t^2-1)+1}{(e^{t^2}-1)^{\frac{3}{2}}}.
\end{align*} 
Proposals are sampled iid from a Rayleigh distribution with the parameter $\sigma=1.3$, which has the density 
\begin{align*}
    g(t)=\frac{t}{\sigma^2}e^{-\frac{t^2}{2\sigma^2}}, 
\end{align*}
for $t\geqslant0$. The probability of accepting a proposed sample from $g$, denoted $x_p$, is 
\begin{align*}
    \frac{f(x_p)}{1.18 \ g(x_p)}. 
\end{align*}
If the proposal is rejected, a new proposal is generated until it is accepted. 
Then, the stochastic representation of Theorem \ref{Th1} is used to obtain a sample from the full distribution.

\subsection{Mat\'ern covariance}
\noindent
The formula for the survival function $E_0$ is given in Table~\ref{Tab:RandE}. One way to simulate is to obtain the inverse of this function $E_0^{-1}$ and simulate $E_0^{-1}(U)$, where $U$ is a uniform random variable on $(0,1)$. 
Since the derivative of $E_0$ is readily available
\begin{align*}
E_0'(t)=-\sqrt{2(\nu-1)}C_\nu t^{\nu-1}
\left(
\frac{t K_{\nu-2}(t)-K_{\nu-1}(t)}{\left(1-C_\nu^2t^{2\nu}K_\nu^2(t)\right)^{1/2}}
+
\frac{C_\nu^2 t^{2\nu+1}K_{\nu-1}(t)K^2_\nu(t)}{\left(1-C_\nu^2t^{2\nu}K_\nu^2(t)\right)^{3/2}}
\right),
\end{align*}
one can seek an importance sampling bound for the density or invert numerically $E_0$ using the Newton–Raphson method. One can also observe that for $\nu=r+1/2$, $r=0,1,2,3,\dots$ the Bessel functions have explicit form
\begin{align*}
K_\nu(u)=\sqrt{\frac\pi{2u}}e^{-u}
\sum_{k=0}^r
\frac{(r+k)!}{(r-k)!k!}(2u)^{-k}
\end{align*}
which significantly reduces the computational cost.
For illustration, consider $\nu=5/2$. Then
\begin{align*}
K_{\nu-2}(u)&=\sqrt{\frac{\pi}{2u}}e^{-u},\\
K_{\nu-1}(u)&=\sqrt{\frac{\pi}{2u}}e^{-u}\left(1+\frac 1 u\right),\\
K_{\nu}(u)&=\sqrt{\frac{\pi}{2u}}e^{-u}
\left(
1+\frac{3}{u}+\frac{3}{u^2}
\right).
\end{align*}
This leads to
\begin{align*}
E_0(t)={\sqrt{3}}\frac{t^2+t}{\sqrt{9e^{2t}- \left(t^2+ 3t +3 \right)^2}}.
\end{align*}
and
\begin{align*}
E_0'(t)
&=E_0(t)\left(\frac{2t+1}{t^2+t}+\frac{(t^2+3t+3)(2t+3)-9e^{2t}}{9e^{2t}-\left(t^2+3 t +3\right)^2}\right).
\end{align*}
Based on this evaluation, the Newton–Raphson method of finding the solution to ${E_0(t)=u}$ has been implemented. A slight modification was needed to ensure the stability of the inversion. It has been used to generate the persistency coefficient approximation presented in Subsection~\ref{subsec:Matern}.


\end{document}